\documentclass[11pt]{article}
  \pdfoutput=1
\usepackage[bookmarks=false]{hyperref}
\usepackage{graphicx,amsmath,amsthm,amssymb,adjustbox,xcolor,multirow,comment}
\numberwithin{equation}{section}
\numberwithin{figure}{section}
\newcommand{\Label}{\label}
\usepackage{fullpage}

\newcommand{\beq}{\begin{equation}}
\newcommand{\eeq}{\end{equation}}

\newcommand{\R}{\mathbb{R}}
\newcommand{\RR}{R}

\newcommand{\ueta}{u} 
\newcommand{\xxi}{x} 
\newcommand{\ttau}{t} 
\newcommand{\ds}{\displaystyle}
 \newcommand{\proofend}{\vrule height 6pt width 6pt depth 1pt}
\def\barr{\begin{array}}
\def\earr{\end{array}}

\newtheorem{thm}{Theorem}[section]

\newtheorem{lem}[thm]{Lemma}
\newtheorem{defn}[thm]{Definition}

\begin{document}
\large
 \centerline{A Scalar Conservation Law  for Plume Migration in Carbon Sequestration}
\ \\
\centerline{Elisabeth Brown and Michael Shearer}https://arxiv.org/submit/1810091/addfiles
%

\thispagestyle{plain}
 
\begin{abstract}  
A quasi-linear hyperbolic partial differential equation with a discontinuous flux models geologic carbon dioxide (CO$_2$) migration and storage \cite{hesse08}. Dual flux curves characterize the model, giving rise to flux discontinuities. One convex flux describes the invasion of the plume into pore space, and the other captures the flow as the plume leaves CO$_2$ bubbles behind, which are then trapped in the pore space. 
We investigate the method of characteristics, the structure of shock and rarefaction waves, and the result of binary wave interactions. The dual flux property introduces unexpected differences between the structure of these solutions and those of a scalar conservation law with a convex flux.  During our analysis, we introduce a new construction of cross-hatch characteristics in regions of the space-time plane where the solution is constant, and there are two characteristic speeds. This construction is used to generalize the notion of the Lax entropy condition for admissible shocks, and is crucial to continuing the propagation of a shock wave if its speed becomes characteristic.  
 \end{abstract}

\section{Introduction}

Some 35.7 billion tonnes of carbon dioxide (CO$_2$) were emitted into the atmosphere in 2014 \cite{jrcpbl}, an increase from the previous year's global CO$_2$ emissions of 32 gigatonnes \cite{huppert14}. In 2000, the Intergovernmental Panel on Climate Change projected a range of estimated emissions from fossil fuel combustion and industrial processes for the year 2020; current emissions are within that annual planning range of 29 to 44 billions tonnes of CO$_2$ \cite{ipcc}. The capture of CO$_2$ before its exodus into the atmosphere seems to be a promising technological solution to reduce the escalating global impact of CO$_2$ emissions. In such a process, gaseous CO$_2$ is collected at industrial sites and power plants, compressed, and injected into geological formations deep underground. Geotechnical evidence suggests that there is a potential subsurface storage capability of 2,000 billion tonnes of CO$_2$ in porous reservoirs worldwide \cite{ipcc}. A goal of future and ongoing carbon dioxide capture and storage projects, such as the Sleipner project  
 beneath the North Sea \cite{sleipner,vella06,zhu15}, is to permanently trap CO$_2$   underground 
  \cite{golding11,nordbotten05}. While a wealth of seismic surveys of the Sleipner project have indicated no signs of leakage \cite{mitccs}, the possibility of escape of the injected CO$_2$ from brine-filled aquifers remains a concern.

During injection, the captured gaseous CO$_2$ is compressed and becomes supercritical; hence, upon release into the porous rock, the sequestered CO$_2$ behaves like a liquid. Since it is less dense than   the ambient brine,   the injected plume rises  within the aquifer \cite{hayek09,hesse08,hesse10}. Appropriate sites for carbon capture and storage projects have an impermeable cap rock in the geological formation that acts as a   barrier to hinder the upward migration of the buoyant plume and keep the CO$_2$ beneath the Earth's surface. Once the plume {rises to} the impermeable upper boundary, the CO$_2$ travels along inclines in the cap's lower surface and spreads through the porous rock as a gravity current. As the plume migrates, it deposits bubbles of CO$_2$  that remain in place. 
The sequestration is successful if all of the CO$_2$ in the plume is deposited  before the plume reaches fractures within  the cap rock that would allow leakage of the plume from the aquifer \cite{golding11,huppert14,silin09,vella06}. 

This mechanism to permanently immobilize CO$_2$ within a porous medium is known as residual trapping.  Capillary forces between the two fluids (brine and supercritical CO$_2$) stably trap bubbles of CO$_2$ within pore spaces. 
Hesse et al.  \cite{hesse08} formulated a quasi-linear hyperbolic partial differential equation with a discontinuous flux to model geologic carbon dioxide migration and storage through residual trapping. 
A striking feature of their   model is that, due to the discontinuous flux,  
the entire CO$_2$ plume is deposited as bubbles in a finite time. 

In this paper, we explore the model in more detail, 
approximating solutions of the Cauchy problem using wave-front tracking. In \S2 we describe the model of \cite{hesse08}, whose key feature is a switch between  two flux functions,  occuring when  the plume changes from  propagating into a region of brine to depositing CO$_2$ droplets. In \S3 we describe novel features of the method of characteristics, and the construction of fundamental solutions of the equation, namely shock waves and rarefaction waves. To establish the admissibility of shock waves, we introduce the notion of cross-hatch characteristics to address the ambiguity of characteristic speeds due to the twin flux functions.  {\S4} includes a detailed description of wave interactions, including some properties that do not occur in conventional scalar conservation laws. In \S5 we construct piecewise constant approximate solutions of the Cauchy problem using expansion shocks in place of rarefaction waves. 
We conclude the paper in \S6 with some remarks.
 
\section{The Two-Flux Model }

In this section, we outline several simplifying assumptions about the aquifer and the nature of the flow, then state the model, a first order partial differential equation with a switch in flux depending on whether, at a given location, the CO$_2$ plume is advancing, or depositing bubbles in its wake. 

\subsection{Model Assumptions}

Subsurface geology often has complicated spatial variability, and three-dimensional models of carbon sequestration require unresolved and difficult issues. To simplify matters, we consider a porous aquifer that is locally uniform in the transverse horizontal direction and analyze the two-dimensional propagation of a cross-section of the flow.  Consider a porous aquifer of constant thickness $H$ beneath an impermeable cap rock sloped at constant angle $\theta$. A buoyant plume of supercritical carbon dioxide, CO$_2$, with height $h(x,t)$ at position $x$ and time $t$ is introduced to the brine-filled aquifer for storage, as shown in Fig.\,\ref{residualbubbles}(a). 
 As in the figure, the CO$_2$ plume is represented by a sharp interface, beglecting  the dissolution of CO$_2$ into the brine   \cite{hesse08,huppert14}. 
The viscosity contrast between the two fluids propels the CO$_2$ plume to invade available pore space as it migrates   as a gravity current \cite{golding11,kestin81,ouyang11}. 

Isolated ganglia of carbon dioxide will be  trapped in a region of the permeable aquifer, with residual surface once the plume recedes,  Fig.\,\ref{residualbubbles}(b). 
Thus, the volume of CO$_2$ within the plume decreases, as the plume migrates and becomes disconnected from the immobilized residual bubbles. It is assumed that pressure within the current is hydrostatic since the advection-dominated migration is mainly horizontal. Within the aquifer, volume is conserved, and the multiphase extension of Darcy's law is applicable in place of conservation of momentum \cite{hayek09,huppert14,juanes10,vella06}. Combining these assumptions with a hyperbolic limit yields a non-dimensional first order partial differential equation given in \cite{hesse08}. 

Let $\ueta=\frac{h}{H}\in\big[\,0\,,1\,\big]$ be the dimensionless height of the CO$_2$ plume, $\ttau$  the non-dimensional advection-dominated time scale, and $\xxi$  the dimensionless spatial variable, based on the initial width $L$ of a  typical plume. The mobility ratio, $\mathcal{M}$, between the supercritical carbon dioxide and the brine depends on permeability and viscosity of each phase; for carbon sequestration, the invading CO$_2$ is more mobile than the ambient brine, so that $\mathcal{M}\geq1$ \cite{hesse08,nordbotten05}. 

The residual surface of immobile CO$_2$ remaining in the wake of the migrating plume   is controlled by a residual trapping parameter, $\varepsilon\in\big[\,0\,,1\,\big]$. Both $\varepsilon$ and $\mathcal{M}$ are constant material properties   \cite{golding11,hesse08,juanes10}. 

\begin{figure}[!h]
\begin{center}
\includegraphics[height=1.65in]{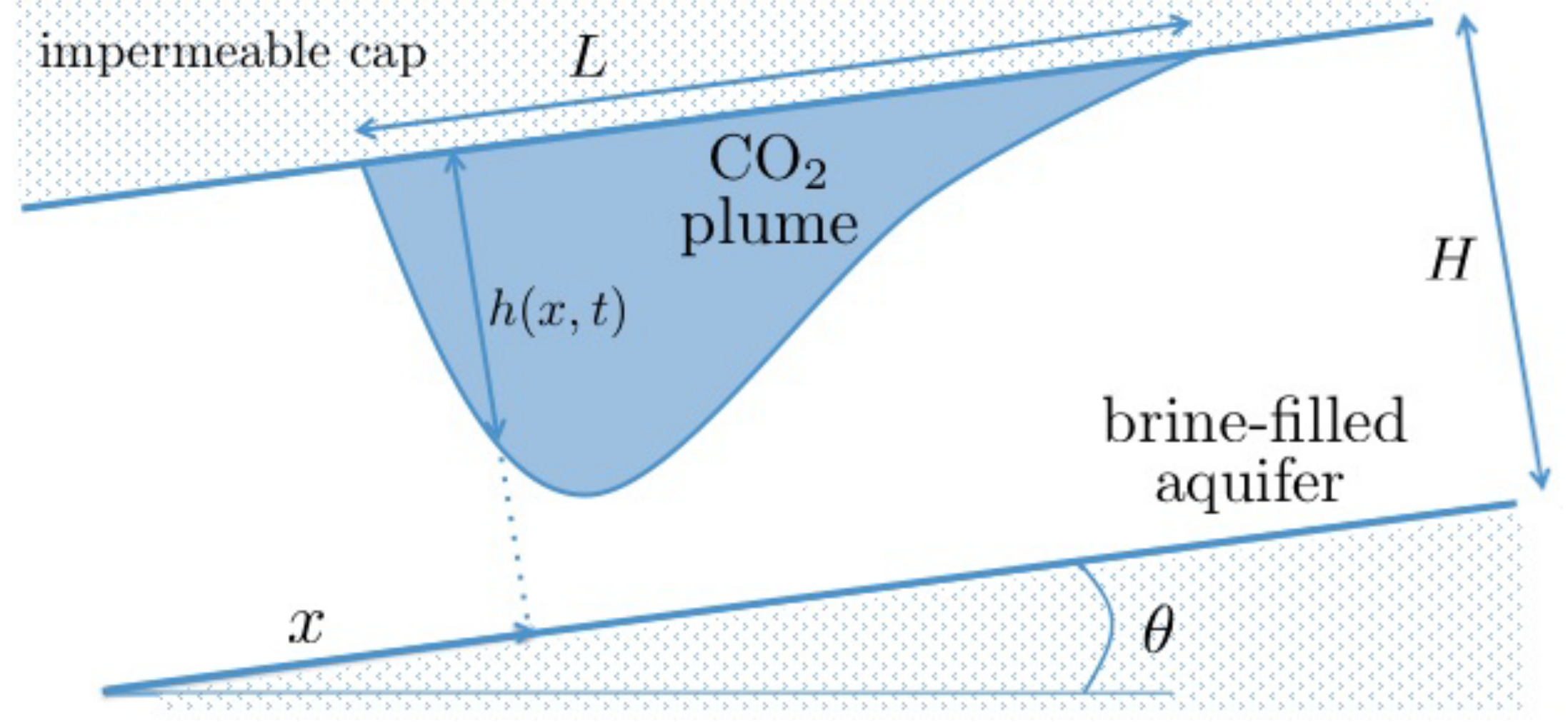}  \includegraphics[height=1.9in]{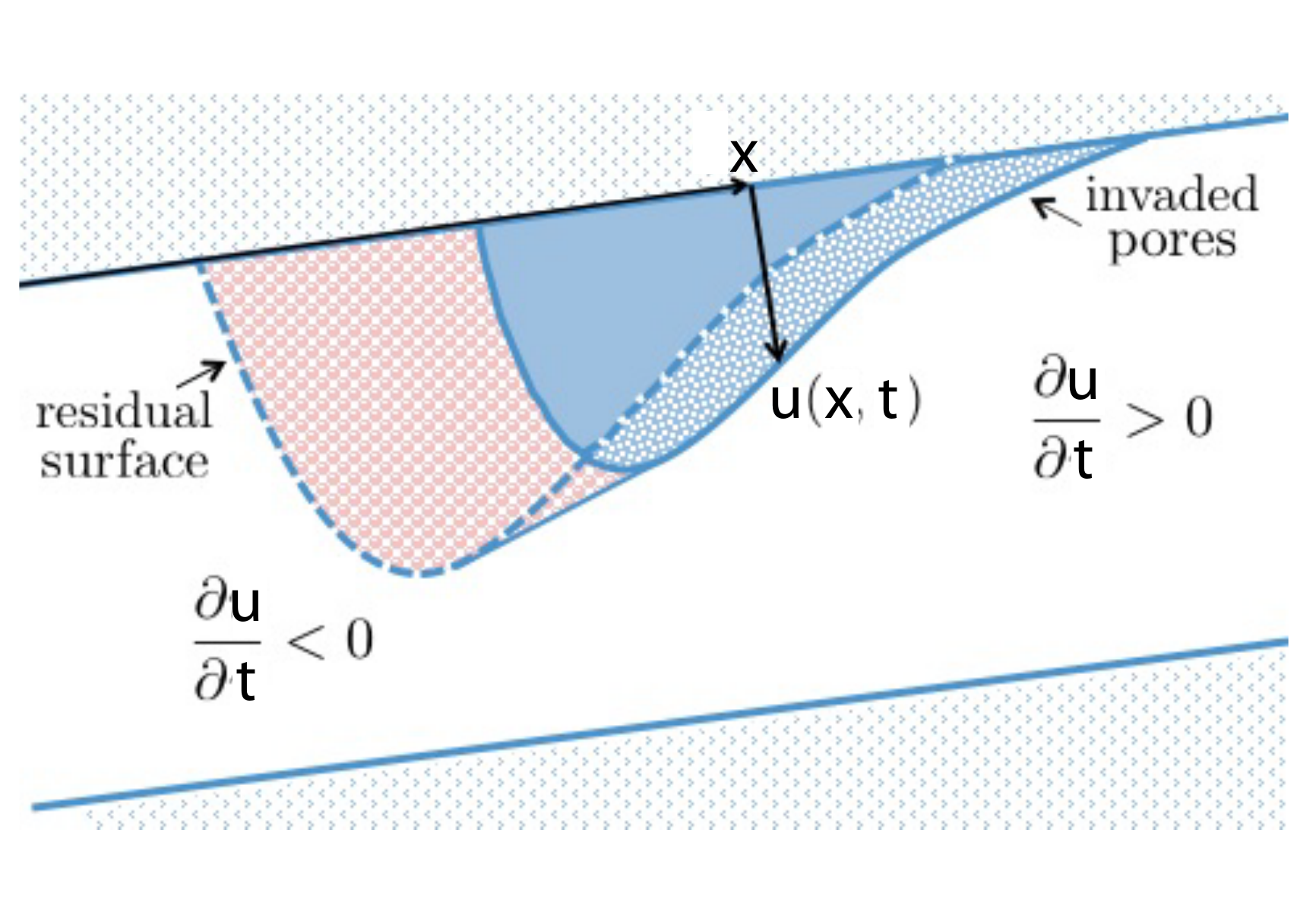} 
(a) \hspace{3in}  (b)
\caption{(a) A CO$_2$ plume in a porous aquifer. (b) Shown in dimensionless variables, a residual area of immobile CO$_2$ remains as the plume migrates to the right.}
\Label{residualbubbles}
\end{center}
\end{figure}

\subsection{Governing Equation}

Hesse, Orr, and Tchelepi \cite{hesse08} modeled the evolution of a gravity current with residual trapping as a scalar equation 
\beq\label{pec1}
\ueta_\ttau + \sigma\, f(\ueta)_\xxi =Pe^{-1}\sigma \big(f(\ueta)\ueta_\xxi\big)_\xxi,
\eeq
in which the  flux $\sigma f$ is a fractional flow rate obtained by eliminating pressure from a version of Darcy's law,
\begin{equation}
\sigma\,f(\ueta)=\,\sigma \, \dfrac{\ueta\,(1-\ueta)}{\;\ueta\,(\mathcal{M}-1)+1\,}\,, 
\Label{flux}
\end{equation} 
and $Pe$ is the Peclet number, representing the ratio of advective and diffusive time scales. 
The    parameter $\sigma\in\big[\,0\,,1\,\big]$ depends on the evolution and is a step function given by
 \begin{equation}
      \sigma =  \left\{
     \begin{array}{ll}
     1-\varepsilon , \ \ &   \mbox{if} \ \ \ueta_\ttau >0\,,\\[2pt]
     1,&    \mbox{if} \ \ \ueta_\ttau <0\,.
     \end{array}
   \right.
   \Label{sig}
\end{equation}
When $\ueta_\ttau>0$, the migrating CO$_2$ is invading new pore spaces, 
whereas when $\ueta_\ttau<0$ the plume is draining, no new trapping locations are sought and the brine invades, isolating bubbles of CO$_2$.

In a sloping aquifer, advection dominates diffusion, and the equation reduces (in the limit $Pe\to \infty$) to the nonlinear conservation law  \begin{equation}
\ueta_\ttau +\sigma\, f(\ueta)_\xxi =0\;, 
\Label{PDE}
\end{equation}
The switch between migration and deposition represented by the parameter $\sigma$  gives rise to discontinuities in the flux. As shown in Fig.\,\ref{ffig}, the lower flux curve describes the invasion of the plume into pore space, and the upper flux captures the flow as the plume leaves CO$_2$ bubbles behind, which are then trapped by brine in the pore space. The characteristic speed is therefore increased during deposition.

Flux functions with discontinuities in space have been previously studied, \cite{chen08,may10,SSShen15}; however,   the flux in this model depends on the sign of $\ueta_\ttau,$ a different kind of discontinuity that introduces new phenomena. 
For $\varepsilon=0$, there is a single flux function; the aquifer has no available pore space to trap CO$_2$, and the plume migrates according  to the classical case in which  the   plume  volume remains fixed and would migrate indefinitely with no deposition.  Typically $\varepsilon\in\big(\,0\,,1\,\big]$ in geologic storage \cite{hesse08,juanes10,nordbotten05,qi09}, and the entire compactly supported plume may be trapped within available pore space after a finite time and within a finite aquifer volume.

\begin{figure}[!h]
\begin{center}
\begin{tabular}{c}
\includegraphics[width=3.5in]{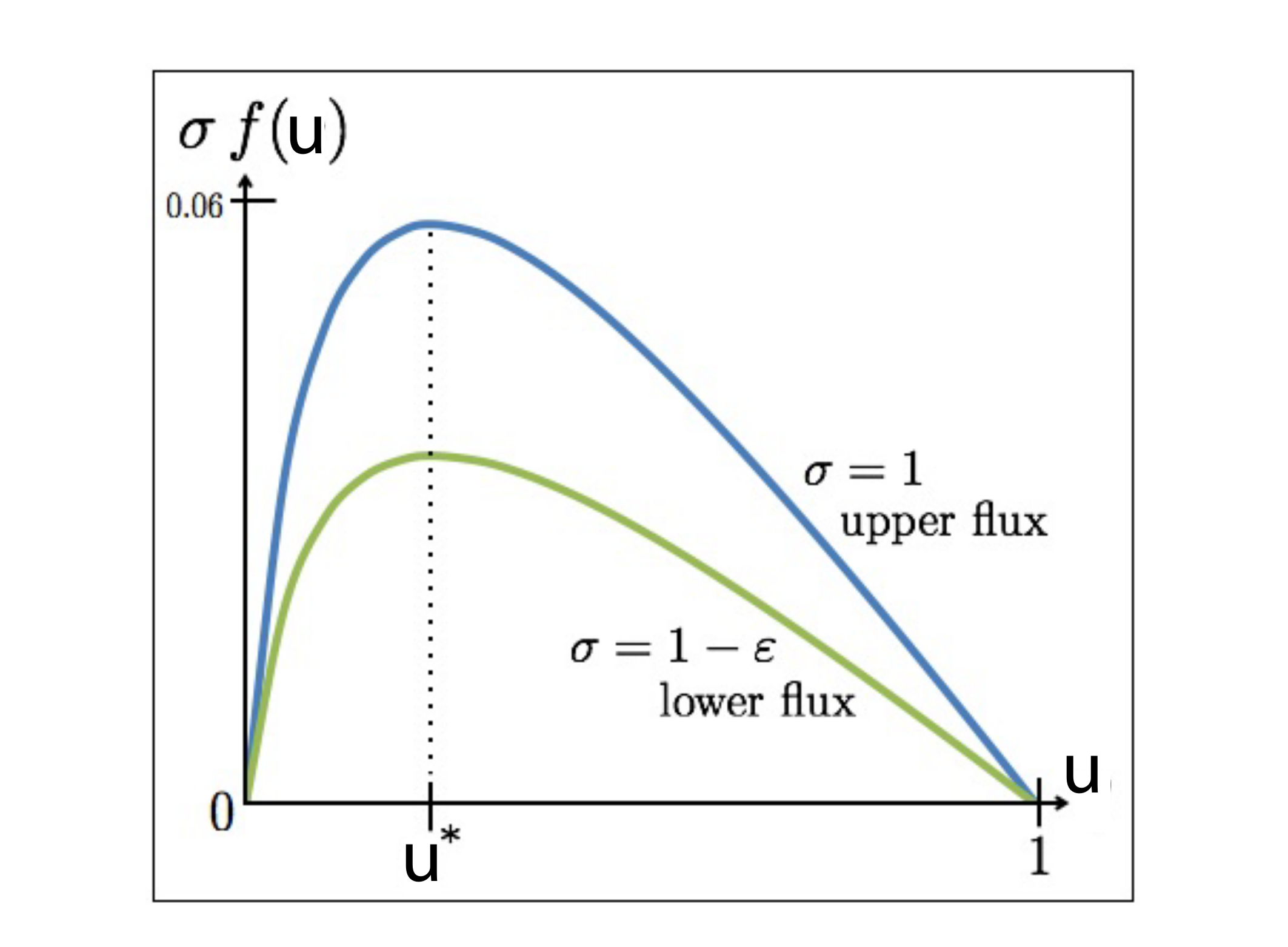}
\end{tabular}
\caption{Dual fluxes \eqref{flux} for $\mathcal{M}=10$ and $\varepsilon=0.4$. Both flux curves attain a maximum value at $\ueta^*= {1}\,\big/\big({1+\sqrt{\mathcal{M}}}\,\big)\,.$  The characteristic speeds satisfy $0<\sigma\,f'(\ueta)<f'(\ueta)$ if $\ueta<\ueta^*$, and   $f'(\ueta)<\sigma\,f'(\ueta)<0$ if $\ueta>\ueta^*$.} 
\Label{ffig}
\end{center}
\end{figure}

\section{Characteristics and Shocks} \Label{charsshocks}

In this section, we consider equation \eqref{PDE} with the switch parameter $\sigma$ given by \eqref{sig}, and assumptions on the flux $f$ consistent with the Hesse et al model \cite{hesse08}:\\

\noindent (H) $f:[0,1]\to \R$ is $C^2, \ f(0)=f(1)=0, \ f''(\ueta)<0, \ 0\leq \ueta\leq 1.$\\

The value $\ueta=\ueta^*$ with $f'(\ueta^*)=0$ plays a significant role in the construction of admissible shock solutions of \eqref{PDE}. For the flux function \eqref{flux}, we have $\ueta^*= {1}\,\big/\big({1+\sqrt{\mathcal{M}}}\,\big)\,.$

We explain the role of the discontinuous switch function $\sigma$ (see \eqref{sig}), in the construction of continuous solutions, shocks and rarefactions. We also resolve an ambiguity, related to the constant regions of $\ueta$ in the characteristic plane, by introducing cross-hatch characteristics. 

\subsection{Method of Characteristics} 
Suppose $\ueta(\xxi,\ttau)$ is a continuous solution of \eqref{PDE} with  initial data $\ueta(\xxi,0)=\ueta_0(\xxi)$ in $C^2.$ For short time, the solution should be obtained by the method of characteristics. However, 
since there are two possible characteristic speeds, $f'(\ueta)$ and $(1-\epsilon)f'(\ueta),$  we have to choose between them, at least in open regions of the $\xxi,\ttau$ plane where $\ueta(x,t)$ is non-constant. Where $\ueta_x(\xxi,t)\neq 0,$ we see that the choice of characteristic speed depends on the sign of $\ueta_x(\xxi,t)$ and the sign of $f'(\ueta),$ since 
\beq\label{etatau}
\ueta_\ttau=-\sigma f'(\ueta)\ueta_\xxi. \  
\eeq
Suppose $\ueta'_0(\bar{\xxi})=0.$ If $\ueta'_0(\xxi)$ is constant in a neighborhood of $\bar{\xxi},$ then the solution is continued to $\ttau>0$ as  that constant, and we introduce cross-hatch characteristics, meaning that both charactistic speeds apply where $\ueta(\xxi,\ttau)$ is constant in an open $\xxi,\ttau$ region. If $\ueta_0$ has an inflection point at $\xxi,$ then the function is either increasing or decreasing at $\xxi,$ and the characteristic speed is uniquely defined. However, if $\ueta_0(\xxi)$ has a maximum or minimum at $\xxi=\bar{\xxi},$ then something interesting happens. Suppose for now that  $f'(u_0(x))>0.$

(1) If $\ueta_0(\xxi)$ has a  minimum at $\xxi=\bar{\xxi},$ with then near $(\xxi,\ttau)=(\bar{\xxi},0),$ the characteristics originating from $\xxi<\bar{\xxi}, \ttau=0$ are slower than those originating from $\xxi>\bar{\xxi}, \ttau=0.$ Consequently, the solution satisfies $\ueta(\xxi,\ttau)=\ueta_0(\bar{\xxi}),$ for $(\xxi,\ttau)$ between the characteristics 
$x=(1-\epsilon)f'(\ueta_0(\bar{\xxi}))t+\bar{\xxi}$ and $x=f'(\ueta_0(\bar{\xxi}))t+\bar{\xxi}.$

(2) If $\ueta_0(\xxi)$ has a  maximum at $\xxi=\bar{\xxi},$ then near $(\xxi,\ttau)=(\bar{\xxi},0),$ the characteristics originating from $\xxi<\bar{\xxi}, \ttau=0$ are faster than those originating from $\xxi>\bar{\xxi}, \ttau=0.$ In this case, the solution has a corner
along a curve $\xxi=\gamma(\ttau),$ with $\gamma(0)=\bar{\xxi}.$ 
\begin{lem} Suppose $\ueta(\xxi,\ttau)$ is a piecewise $C^2$ solution of equation \eqref{PDE} satisfying $\ueta(\xxi,0)=\ueta_0(\xxi),$ where $\ueta_0\in C^2(\R).$
If $\ueta_0(\xxi)$ has a maximum at $\xxi=\bar{\xxi},$  and $f'(\ueta_0(\bar{\xxi}))>0,$ the maximum propagates as a corner $\xxi=\gamma(t)$ in the graph of $\ueta(\xxi,\ttau), \ \ttau>0,$ satisfying
\beq\label{gammaprime}
\gamma'(\ttau)=f'(\ueta(\gamma(\ttau),\ttau))\left(1+\epsilon\frac{\ueta_x^+(\ttau)}{\ueta_x^-(\ttau)-\ueta_x^+(\ttau)}\right), \ t>0,
\eeq
where   $\ueta_x^\pm(\ttau)=\ueta_\xxi(\gamma(\ttau)^\pm,\ttau);$ 
$  \gamma(0)=\bar{\xxi}, \ \gamma'(0)= (1-\frac\epsilon 2)c, \ c=f'(\ueta_0(\bar{\xxi})).$
\end{lem}
\noindent{\bf Proof:} \ 
To derive an ODE for $\gamma(\ttau),$ we differentiate 
the continuity condition
$$
\ueta(\gamma(\ttau)^-,\ttau)=\ueta(\gamma(\ttau)^+,\ttau)
$$
with respect to $\ttau,$ and use the identity \eqref{etatau}.  After some manipulation, we establish \eqref{gammaprime} for $\ttau>0,$ where necessarily $\ueta_x^-(\ttau)>0>\ueta_x^+(\ttau).$ Note that away from $x=\gamma(t),$ the solution is determined from the method of characteristics. Thus, $u_x^\pm(t) $ depend implicitly on $\gamma(t):$
$$
u_x^\pm(t) = \frac{u_0'(\gamma-\sigma f'(u(\gamma,t))}{1+\sigma u_0'( \gamma-\sigma f'(u(\gamma,t))f''(u_0(\gamma-\sigma f'(u(\gamma,t)))t},
$$
where $\gamma=\gamma(t)$ and $\sigma=1-\epsilon, 1$ for $u_x^\pm$ respectively. 

However,   equation  \eqref{gammaprime}  has a singular limit  as $\ttau\to 0,$ since $  \ueta_x^\pm(0) =0.$
Let $c=f'(\ueta_0(\bar{\xxi})),$ and note that $f'(\ueta(\gamma(\ttau),\ttau))=c$ to leading order as $\ttau\to 0.$  Without loss of generality, we assume that $c>0.$
Similarly, since $\ueta_0(\xxi)$ is $C^2$ and has a maximum at $\xxi=\bar{\xxi},$ if $\ueta_0(\xxi_L)=\ueta_0(\xxi_R),$ with 
$\xxi_L<\bar{\xxi}<\xxi_R,$ then to leading order, $\xxi_L-\bar{\xxi}=\bar{\xxi}-\xxi_R.$ Now consider the solution $\ueta(\xxi,\ttau).$ It is determined at the maximum from two different characteristics, that meet at $\xxi=\gamma(\ttau).$
If the two characteristics emanate from $\xxi_L<\xxi_R,$ then $\gamma=f'(\ueta)\ttau+\xxi_L=(1-\epsilon)f'(\ueta)\ttau+\xxi_R.$ Thus, to leading order near $\ttau=0,\ \gamma =c\ttau+2\bar{\xxi}-\xxi_R =(1-\epsilon)c\ttau +\xxi_R.$
Solving the second equation, we have $\xxi_R=\frac12 c \ttau\epsilon+\bar{\xxi}.$ Hence, $\gamma=\bar{\xxi}+(1-\frac\epsilon 2)c\ttau$ to leading order. Thus, as $\ttau\to 0, \ \gamma'(\ttau)\to  (1-\frac\epsilon 2)c.$ That is, the initial speed of the corner, at the local maximum of $\ueta(\xxi,\ttau)$ (with respect to $\xxi$) is the average of the two characteristic speeds $c$ and $(1-\epsilon)c.$  \ \proofend

\noindent{\bf Remarks} \ 1. \ If $c=f'(\ueta_0(\bar{\xxi}))<0,$ a corresponding argument applies, but the propagation is to the left. In this case, we have 
$$
\gamma'(\ttau)=f'(\ueta(\gamma(\ttau),\ttau))\left(1+\epsilon\frac{\ueta_x^-(\ttau)}{\ueta_x^-(\ttau)-\ueta_x^+(\ttau)}\right), \ t>0.
$$

2. The functions $u^\pm_x(t) $ depend implicitly on $\gamma$ as follows.For $f'(u_0(\bar{x}))>0,$ we have (for $x$ near $\bar{x}$)
$$
u(x,t)=\left\{\barr{ll}
u_0(x-f'(u)t), \quad& x\leq \gamma(t)\\[6pt]
u_0(x-rf'(u)t), \quad& x \geq \gamma(t),
\earr
\right.
$$
and note that $u(x,t)$ is continuous, at least over some finite time interval $0\leq t\leq T.$ 
Differentiating with respect to $x,$ we have $u_x^+(t)=u_0'/(1+ru_0'f''(u)t),$ where 
$u_0'=u_0'(\gamma(t)-rf'(u)t),$ and $u=u(\gamma(t),t),$ and a similar expression for $u_x^-(t),$ but dropping $r$ from both expressions.  

\subsection{Cross-hatch Characteristics} \Label{crosshatchsec}

Since the switch parameter $\sigma$ is not defined when $\ueta_\ttau=0$, the characteristic speed is not well defined in regions of the characteristic plane where the solution is constant. To resolve this, we include 
characteristics determined by both flux curves at each point where $\ueta_\ttau=0$\,; we refer to them as \textit{cross-hatch characteristics} since they form a cross-hatch pattern in regions where $\ueta$ is constant (see Fig.\,\ref{shocking}(a)). 

The two possible characteristic speeds are $\sigma f'(\ueta), \ \sigma=1$ or $\sigma=1-\varepsilon.$ We refer to the larger or greater characteristic speed as the {\it faster} speed, and the other characteristic speed as the {\it slower} speed. {\it Thus, the faster speed is $f'(\ueta)$ (and hence positive) if and only if $\ueta<\ueta^*.$ } To clarify further, when $\ueta>\ueta^*,$ we have  $ f'(\ueta)<0,$ so that the faster  speed is $(1-\epsilon)f'(\ueta)$ since it is greater than  $f'(\ueta)$ in this case.

\subsection{Shocks} \Label{shockssection}
The definition of weak solution for equation \eqref{PDE} does not follow the usual pattern of multiplication by a test function and integration by parts. To see that the usual procedure is problematic, we rewrite the equation as
\beq\label{PDEa}
\ueta_\ttau +\sigma(u_t) f(\ueta)_\xxi =0, \quad \sigma(u_t)=1-\epsilon H(u_t),
\eeq
where $H$ is the Heaviside step function. This form highlights the difficulty of interpreting the equation in the sense of distributions, as both $\sigma(u_t)$ and $f(u)_x$ may be singular.
However, if $u(x,t)$ has only jump discontinuities, then although $u_t$ is singular, namely a delta function, the definition of $H(u_t)$ can be extended by $H(a\delta(x))=1$ if $a>0,$ and $H(a\delta(x))=0$ if $a\leq 0.$ In this way, the notion of solution can be extended to piecewise smooth functions. 

 To define piecewise smooth solutions with jump discontinuities, it is enough to consider a piecewise constant jump discontinuity 
\begin{equation}
       \ueta\left(\xxi,\ttau\right)= \left\{
     \begin{array}{ll}
     \ueta_L, &\quad\xxi<\Lambda \ttau,\\[2pt]
     \ueta_R,&\quad \xxi>\Lambda \ttau
     \end{array}
     \right.\Label{shocksoln}
\end{equation}
propagating with speed $\Lambda.$ 
 Since $\sigma$ in \eqref{PDEa} is selected by the sign of $\ueta_\ttau$, we set $\sigma=1$ if $\ueta$ jumps down across the shock as time increases; otherwise, if the jump is up, we set $\sigma=1-\varepsilon$. This fixes the value of $\sigma$, and we can write the Rankine-Hugoniot jump condition, 
\begin{equation}
\Lambda \,=\, \dfrac{\;\sigma \,\big[\,f(\ueta_R)-f(\ueta_L)\,\big]\;}{\ueta_R-\ueta_L} 
\Label{shockspeed}
\end{equation}
Hesse et al. \cite{hesse08} justified the choice of $\sigma$ in a slightly different way by including dissipative terms (in \eqref{pec1}) that smooth the shock.

For a scalar conservation law with a single flux function, admissible shocks satisfy the Lax entropy condition, requiring characteristics to enter the shock on both sides  \cite{lax}.  Here, with two fluxes, we specify shock admissibility as follows:

{\defn The shock wave \eqref{shocksoln} is admissible if and only if  the faster  characteristics  enter the shock from both sides.}

\ \\
{We argue that   \eqref{shocksoln} is an admissible shock if and only if  $\ueta_L<\ueta_R,$ 
 just as it would be for a scalar conservation law with a convex flux. 
As shown in Fig.\,\ref{shocking}(a), $\ueta_\ttau<0$ across an admissible forward shock (i.e., with $\Lambda>0$), so   that   $\sigma=1.$ Consequently, not only is  $\Lambda$ determined from the upper flux curve, but also the faster characteristics  enter the shock,  see Fig.\,\ref{shocking}(b). 
For an admissible backward shock, with  $\Lambda<0$, we have $\sigma=1-\varepsilon,$ and the shock is admissible if and only if the characteristics found on the lower flux curve impinge on the shock on the right, because they are the less negative characteristics, and enter the shock on the left because either they are the fast characteristics (if $\ueta_L>\ueta^*$), or both families have positive speed  (if $\ueta_L<\ueta^*$), as shown in Fig.\,\ref{backward}.   Once again, this amounts to the condition $\ueta_L<\ueta_R,$ but there is an important point regarding the slower characteristics, which necessarily enter the shock on the right, but may leave on the left. 

\begin{lem}\label{lemma1}
The only characteristics that can leave an admissible shock belong to the slower family, and are on the left of the shock. 
\end{lem}

\textit{Proof}: Consider an admissible shock \eqref{shocksoln}.   If $\ueta_R<\ueta^*$, the faster characteristic speed is on the upper flux, so $f'(\ueta_R)<\Lambda$ is required for admissibility.  Thus, $(1-\varepsilon)\,f'(\ueta_R)<\Lambda$ also. Hence, both characteristics on the right impinge on the shock. If $\ueta_R > \ueta^*$, the faster characteristic speed  is on the lower flux curve, so admissibility requires $(1-\varepsilon)\,f'(\ueta_R)<\Lambda$. Since $\ueta_R>\ueta^*$, $f'(\ueta_R)<(1-\varepsilon)\,f'(\ueta_R)$, and the slower characteristic on the right also enters the shock. Hence, both characteristics on the right always impinge on an admissible shock.
Since the faster characteristics are required to enter the shock on the left, only the slower characteristics on the left can leave the shock. \proofend

\vskip11pt
It is perhaps instructive to understand when the slower characteristics leave an admissible shock. 
If $\ueta_L<\ueta^* $   in a backward admissible shock, both characteristics on the left have positive speed but the shock speed is negative, so both characteristics on the left must enter the shock.

When $\ueta_L<\ueta^*$ in a forward admissible shock, the faster characteristic entering the shock from the left has speed $f'(\ueta_L)>\Lambda\,,$ since $\sigma=1$ in \eqref{shockspeed} for a forward shock.   If $\Lambda<(1-\varepsilon)\,f'(\ueta_L)$, the slower characteristics  will also impinge on the forward shock; however, it is possible that $(1-\varepsilon)\,f'(\ueta_L)<\Lambda$\,, in which case the slower characteristics on the left emanate from the shock. Similarly, if $\ueta^*<\ueta_L$, an admissible   shock requires $0>(1-\varepsilon)\,f'(\ueta_L)>\Lambda$ since $\sigma=1-\varepsilon$ in \eqref{shockspeed}. The more negative characteristic speed $f'(\ueta_L)$  may or may not satisfy $f'(\ueta_L)>\Lambda$\,, so  the slower characteristics on the left can leave the shock. 

In summary, since the faster characteristics must impinge on the shock from both sides, the slower characteristics on the right also enter the shock, but the slower characteristics on the left can leave the shock. Fig.\,\ref{catalogclean2d}(b) illustrates the latter behavior of the characteristics. 
\vskip11pt



\begin{figure}[!h]
\begin{center}
\begin{tabular}{cc}
\includegraphics[height=2.2in]{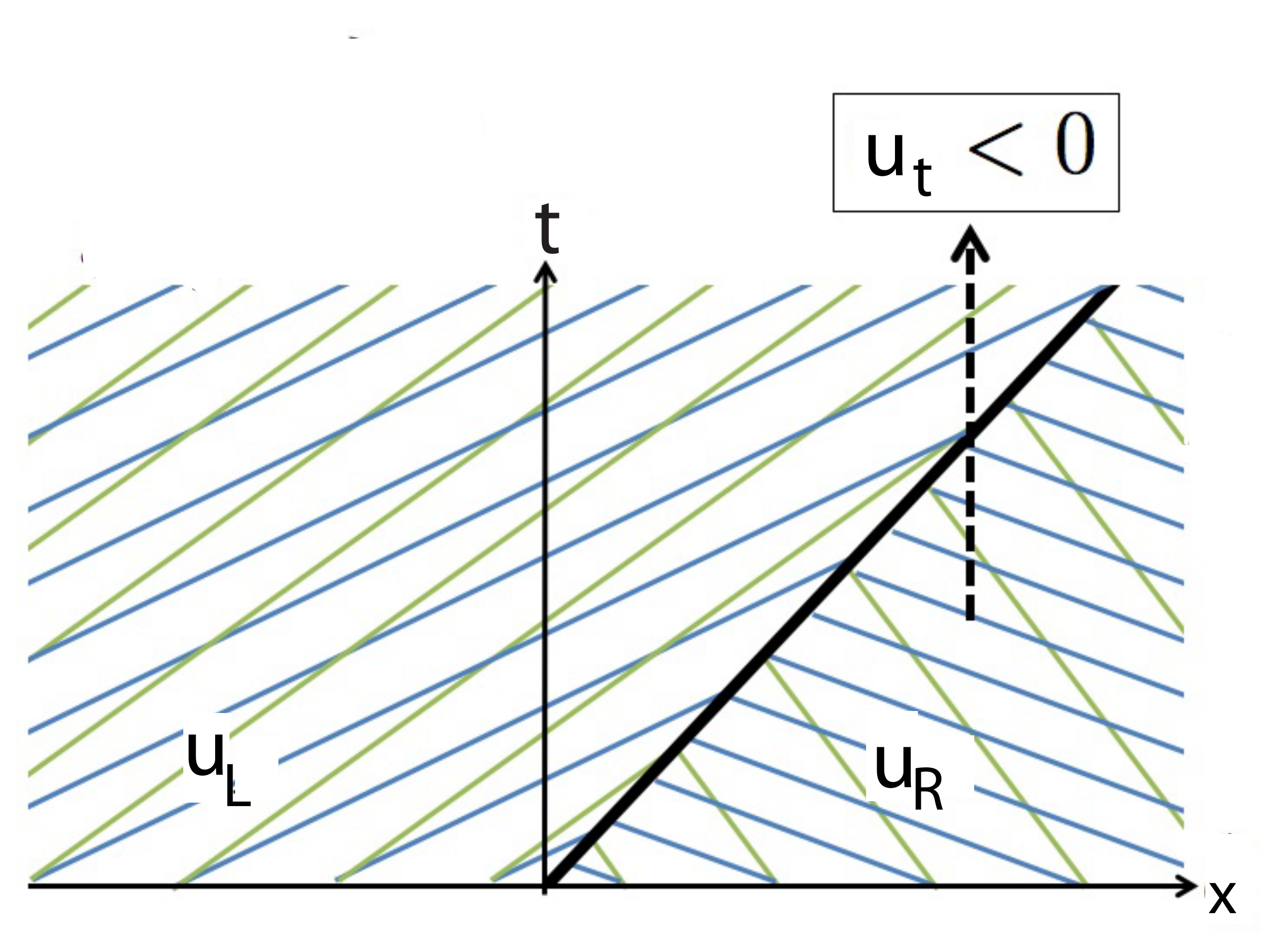}  & \includegraphics[height=2.2in]{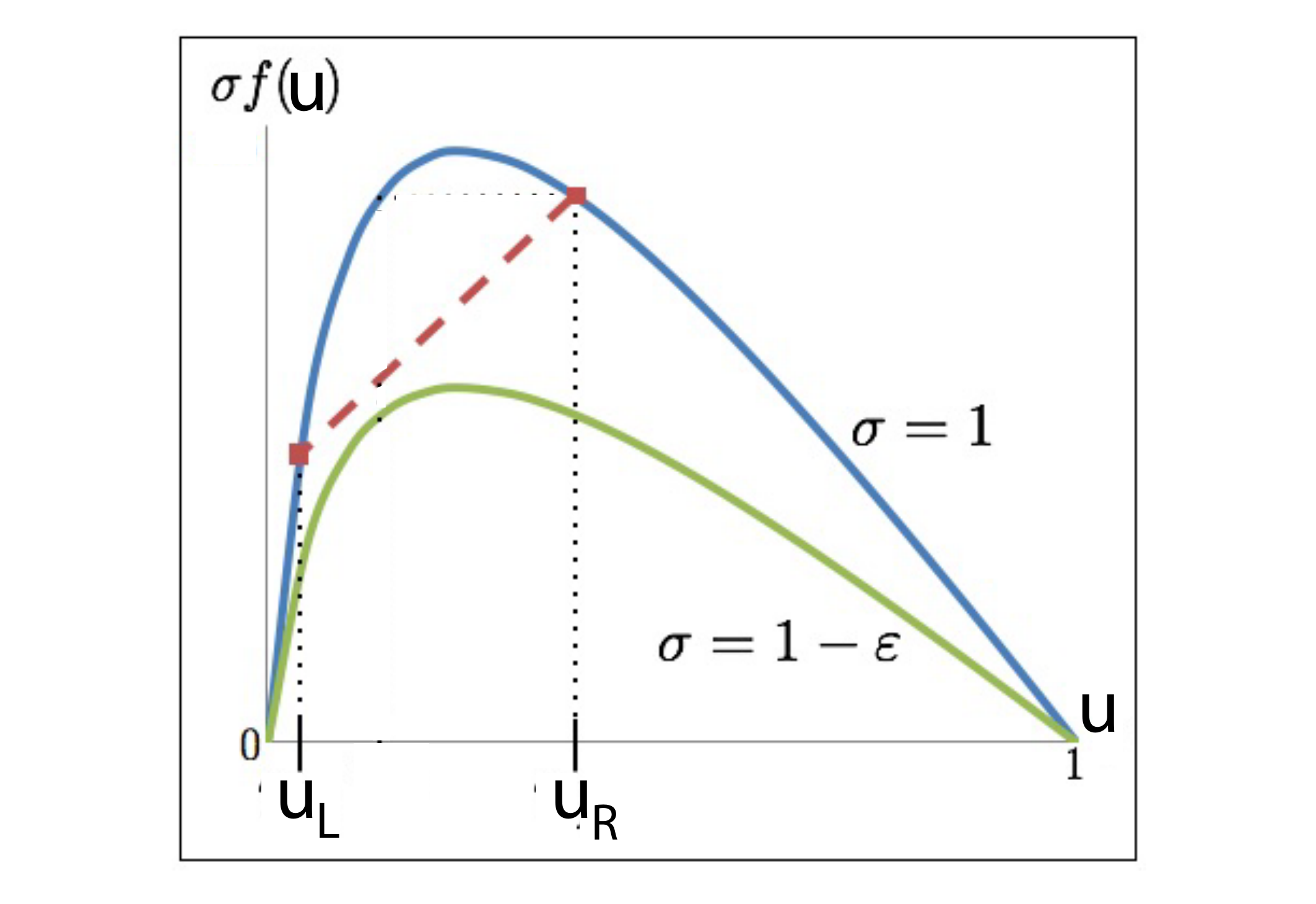} \\
(a) & (b)
\end{tabular}
\caption{A forward shock with $\Lambda>0.$ (a) Characteristic plane with cross-hatch characteristics in constant regions. (b) Shock speed determined from the upper flux curve.}
\Label{shocking}
\end{center}
\end{figure}

 \begin{figure}[!h]
\begin{center}
\includegraphics[height=2.2in]{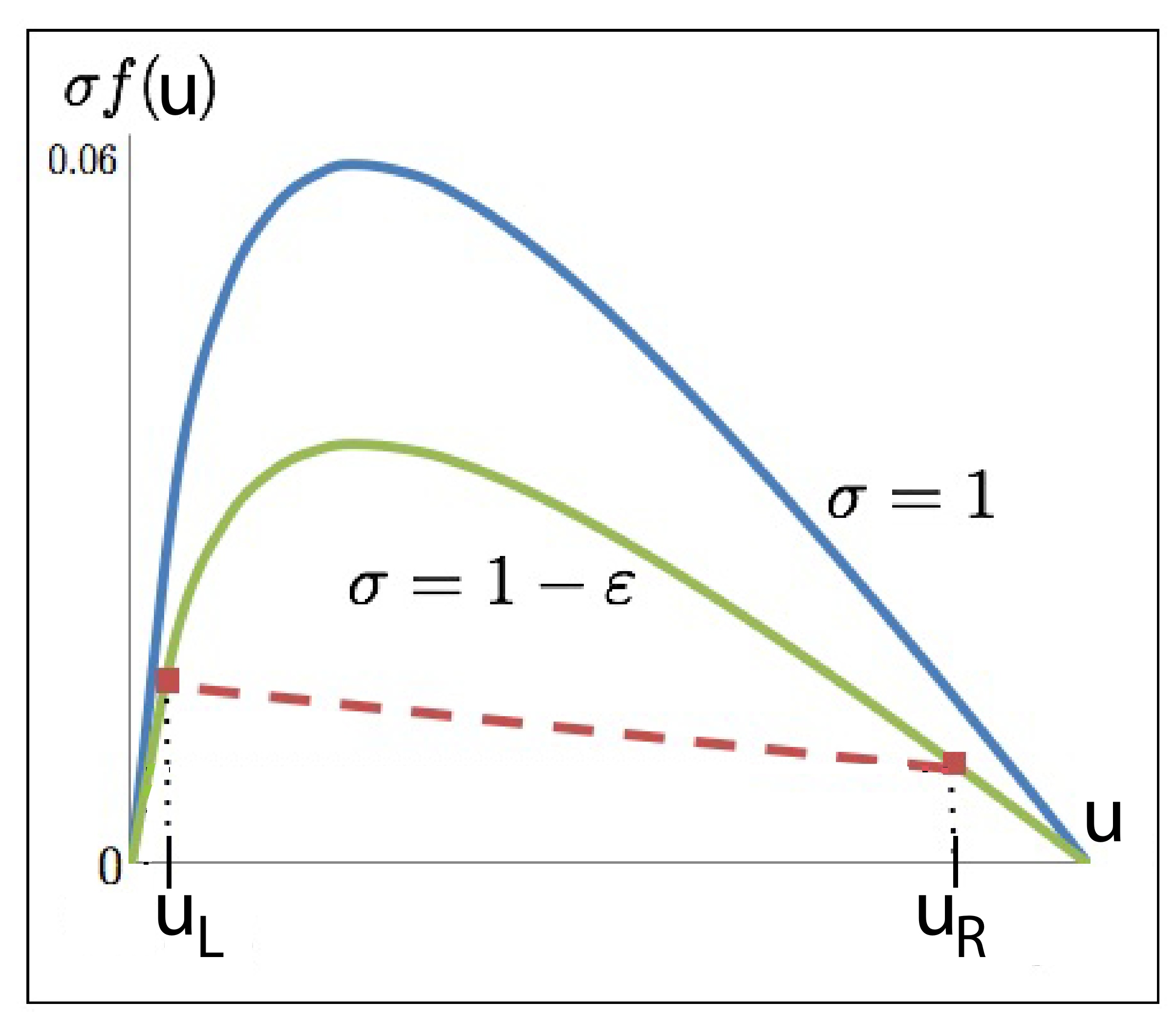} 
\caption{Backward shock, $\Lambda<0\,,$ for which $\ueta_\ttau>0$\,, so the shock speed  is found from the lower flux curve. }
\Label{backward}
\end{center}
\end{figure}

\subsection{Expansion Shocks} 

Expansion shocks are shock wave solutions of \eqref{PDE} that are inadmissible. We characterize them here because we will need them  {in \S5} as approximations to rarefactions in wave-front tracking. For scalar conservation laws, expansion shocks have characteristics leaving in forward time on both sides. 
Here we define a discontinuous function \eqref{shocksoln} to be an expansion shock if  it satisfies the Rankine-Hugoniot jump condition \eqref{shockspeed} and the slower characteristics on each side emanate from the shock.  The latter condition  is equivalent to 
$\ueta_R<\ueta_L.$
Then the faster characteristics on the right also leave the shock, but the faster characteristics on the left may or may not enter the shock, as shown in Fig.\,\ref{expansion}. 

\begin{figure}[!h]
\begin{center}
\begin{tabular}{cc} 
\includegraphics[height=2.in]{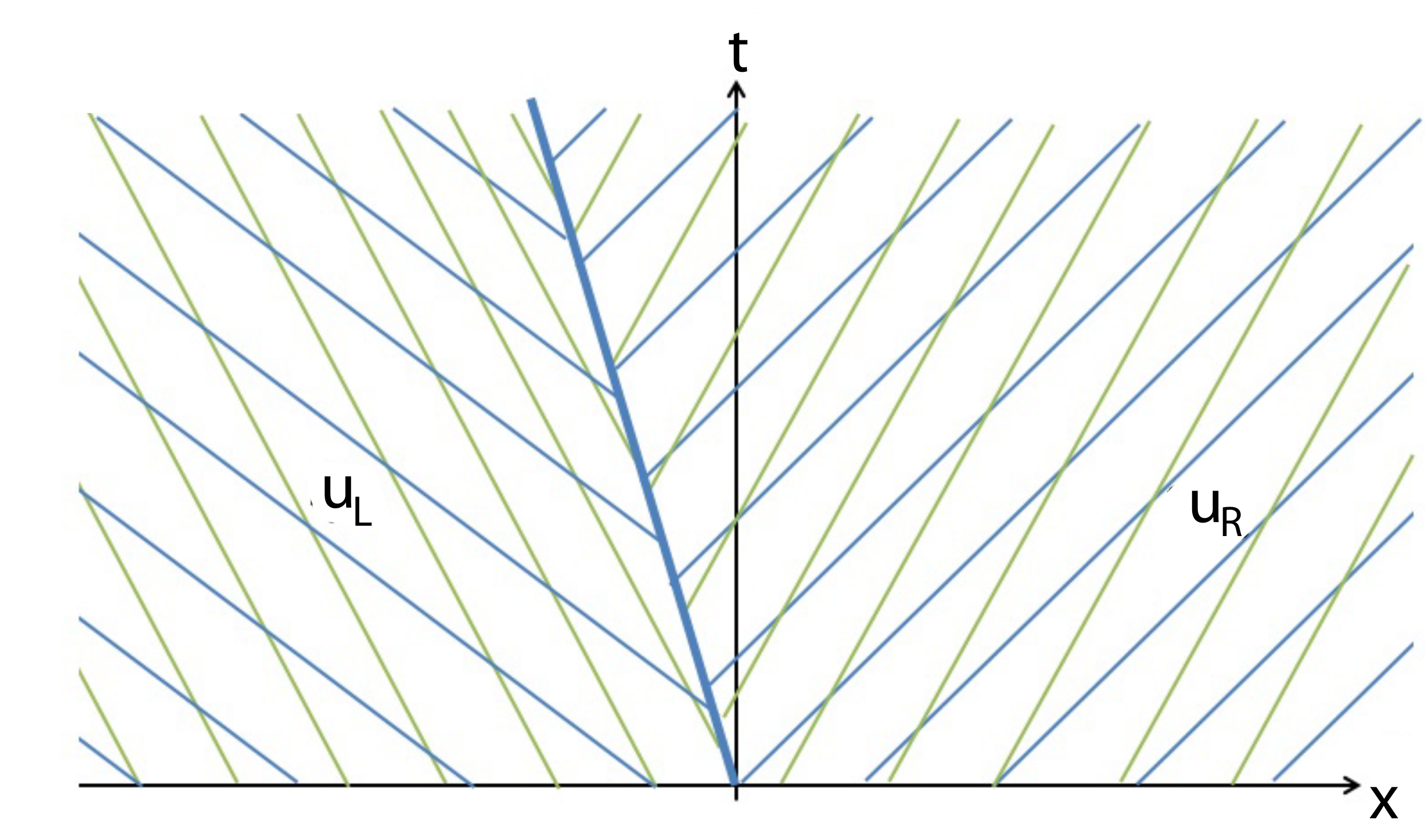}\; & \includegraphics[height=2.in]{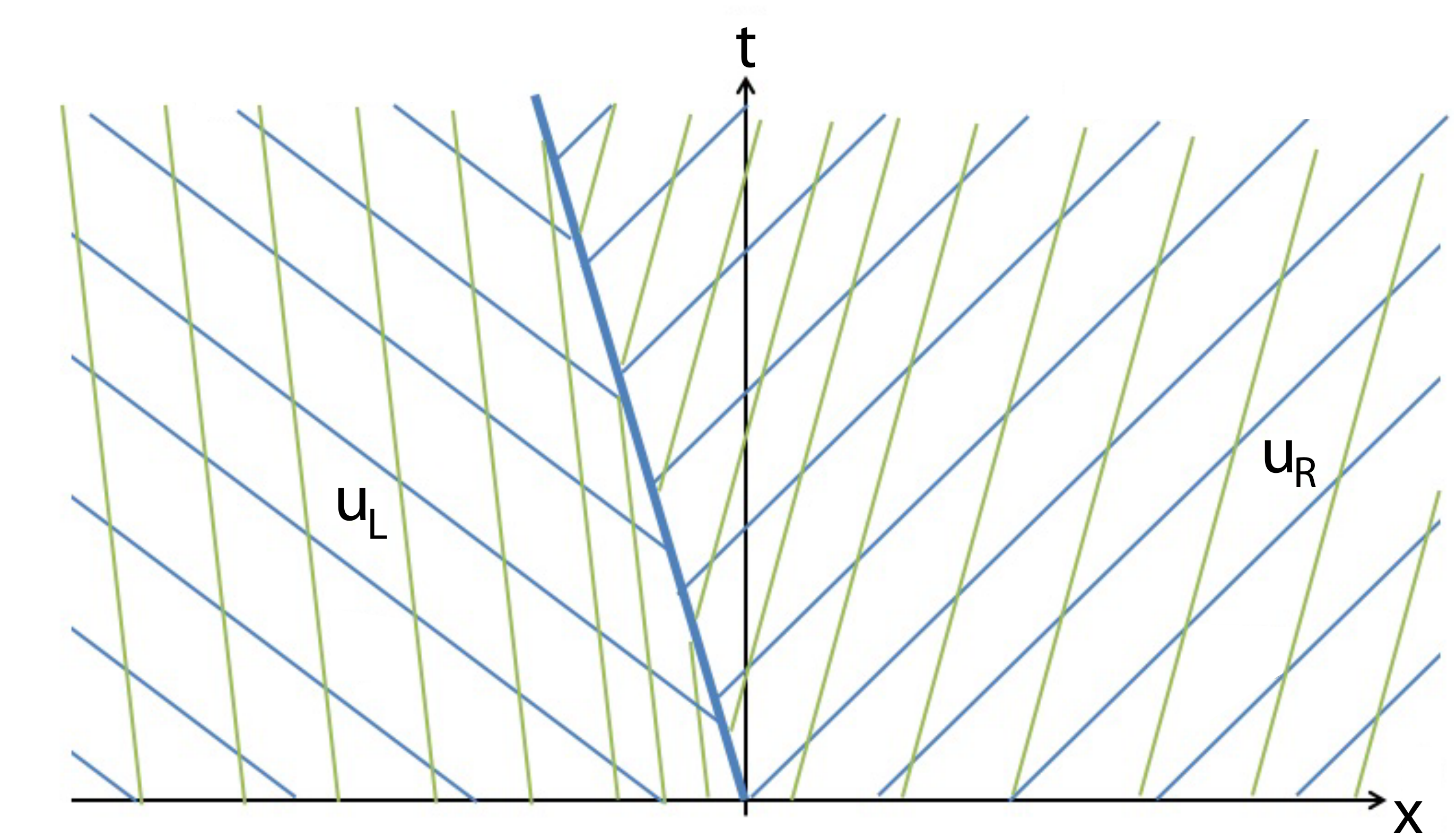} \\ 
(a) \,& (b)\;\; 
\end{tabular}
\caption{Backward expansion shock. (a) All characteristics leave the shock. (b) The faster characteristics on the left enter the shock.}
\Label{expansion}
\end{center}
\end{figure}

\subsection{Rarefactions} \Label{raresection}

Centered rarefaction fans are continuous weak solutions of \eqref{PDE} obtained via the method of characteristics and have the form
\begin{equation}
       \ueta\left(\xxi,\ttau\right)= \left\{
     \begin{array}{lc}
     \ueta_L\,, &\quad\dfrac{\xxi}{\ttau}<\sigma\,f'(\ueta_L)\,,\\[6pt]
     \check{\ueta}\left(\dfrac{\xxi}{\ttau}\right), &\quad \sigma\,f'(\ueta_L)\leq\dfrac{\xxi}{\ttau}\leq\sigma\,f'(\ueta_R)\,,\\[6pt]
     \ueta_R\,,&\quad \sigma\,f'(\ueta_R)<\dfrac{\xxi}{\ttau}\,,
     \end{array}
     \right.\Label{raresoln}
\end{equation}
in which, 
the function $\check{\ueta}$ is given implicitly by $ y=\sigma f'(\check{\ueta}(y)).$


  The rarefaction in Fig.\,\ref{exprare}(a) has both forward and backward characteristics {with speeds that} depend on the value of $\sigma\,,$ as explained in the figure caption.  In Fig.\,\ref{exprare}{(b)} we show a rarefaction wave approximated by three expansion shocks; from left to right, the expansion shocks have increasing speeds. In Fig.\,\ref{rare}(a), we show the construction of the rarefaction wave, resolving the initial  step down in $\ueta\,,$ using both flux functions, and in Fig.\,\ref{rare}(b) we show the corresponding plume profile.
\begin{figure}[!h]
\begin{center}
\begin{tabular}{cc}
\includegraphics[width=3in]{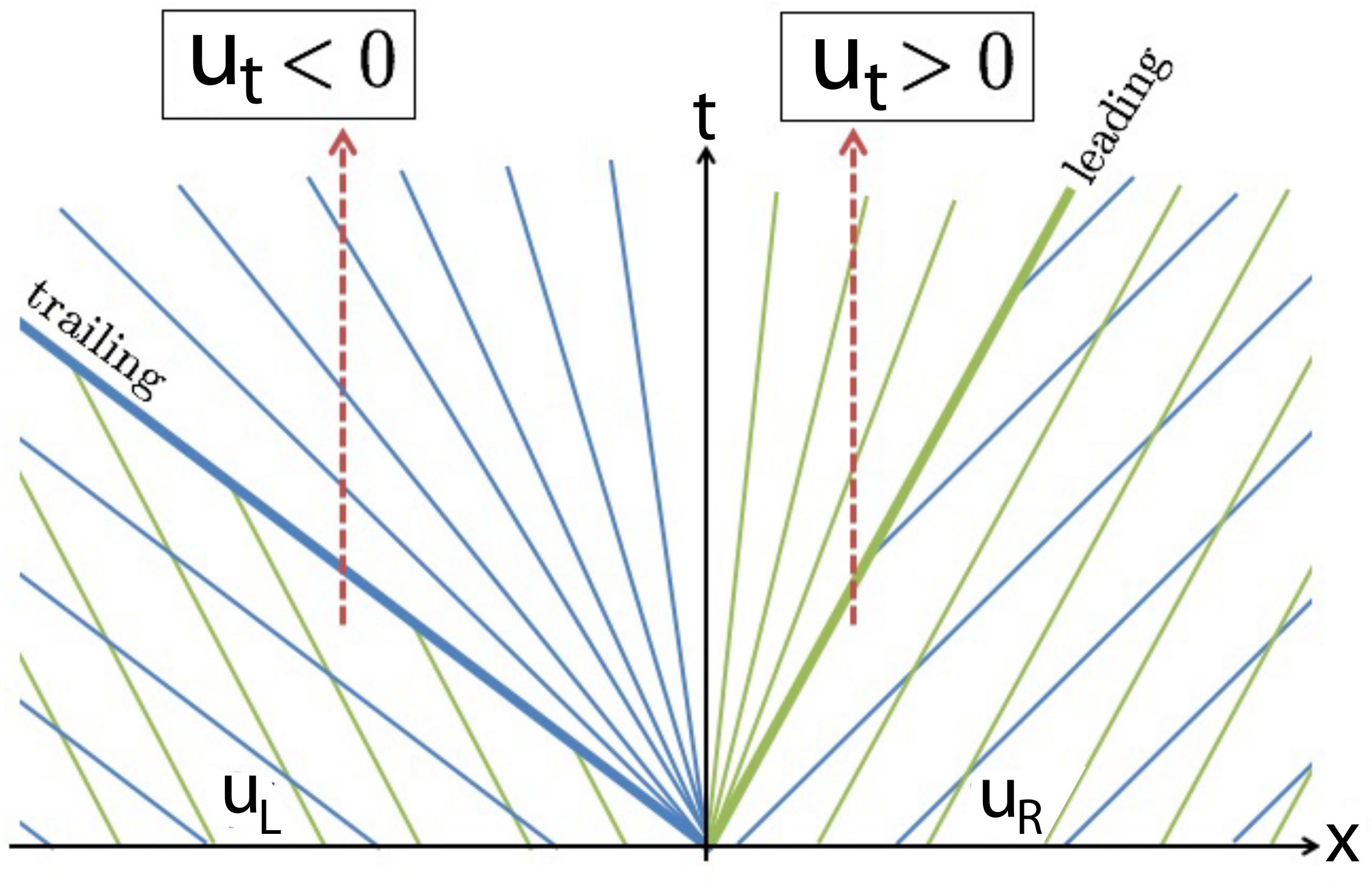} & \includegraphics[width=3in]{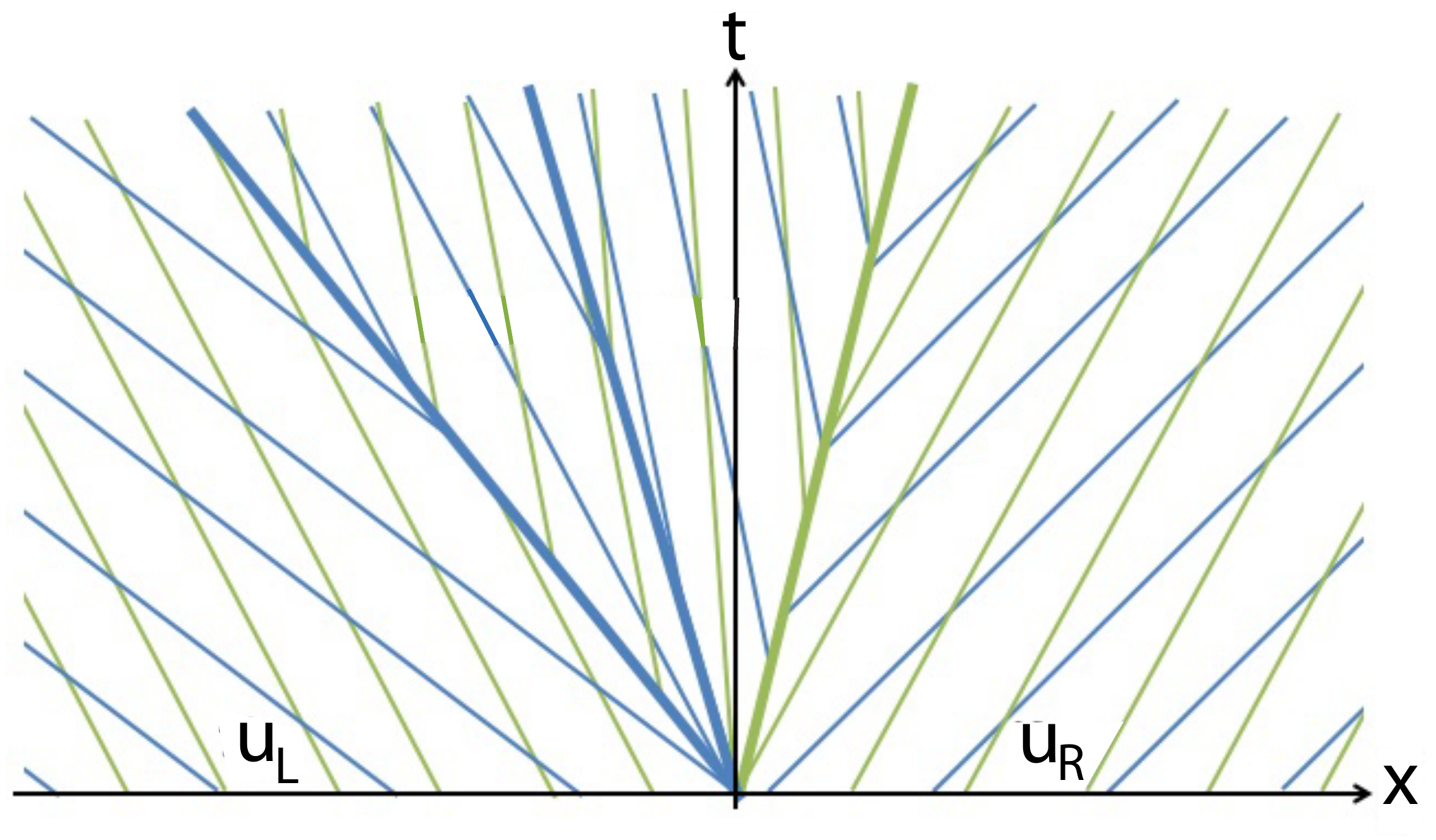} \\
\,\;\;(a) & \;\;(b) 
\end{tabular}
\caption{(a) Rarefaction wave with $\ueta_R<\ueta^*<\ueta_L$. Since $\ueta$ is necessarily  
{decreasing} 
 from left to right in the rarefaction wave, we have that $\ueta_\ttau<0$  left of the $\ttau$ axis,  so that $\sigma=1$. To the right, $\ueta_\ttau>0,$ so that $\sigma=1-\varepsilon$ there. (b) Three expansion shocks approximating the rarefaction wave of (a).   }
\Label{exprare}
\end{center}
\end{figure}


\begin{figure}[!h]
\begin{tabular}{l}
\includegraphics[height=2in]{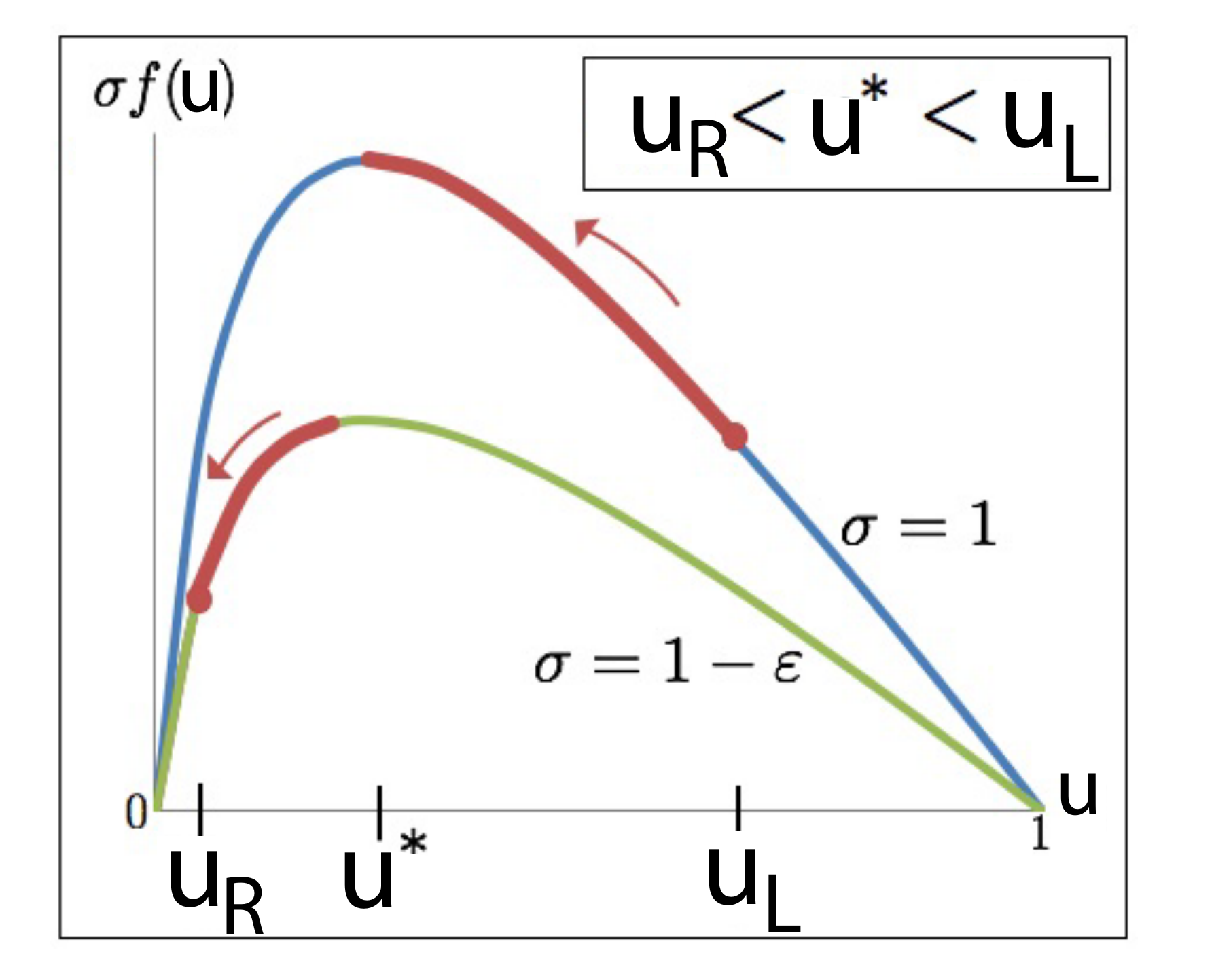} \hspace{.2in} \includegraphics[height=2.3in]{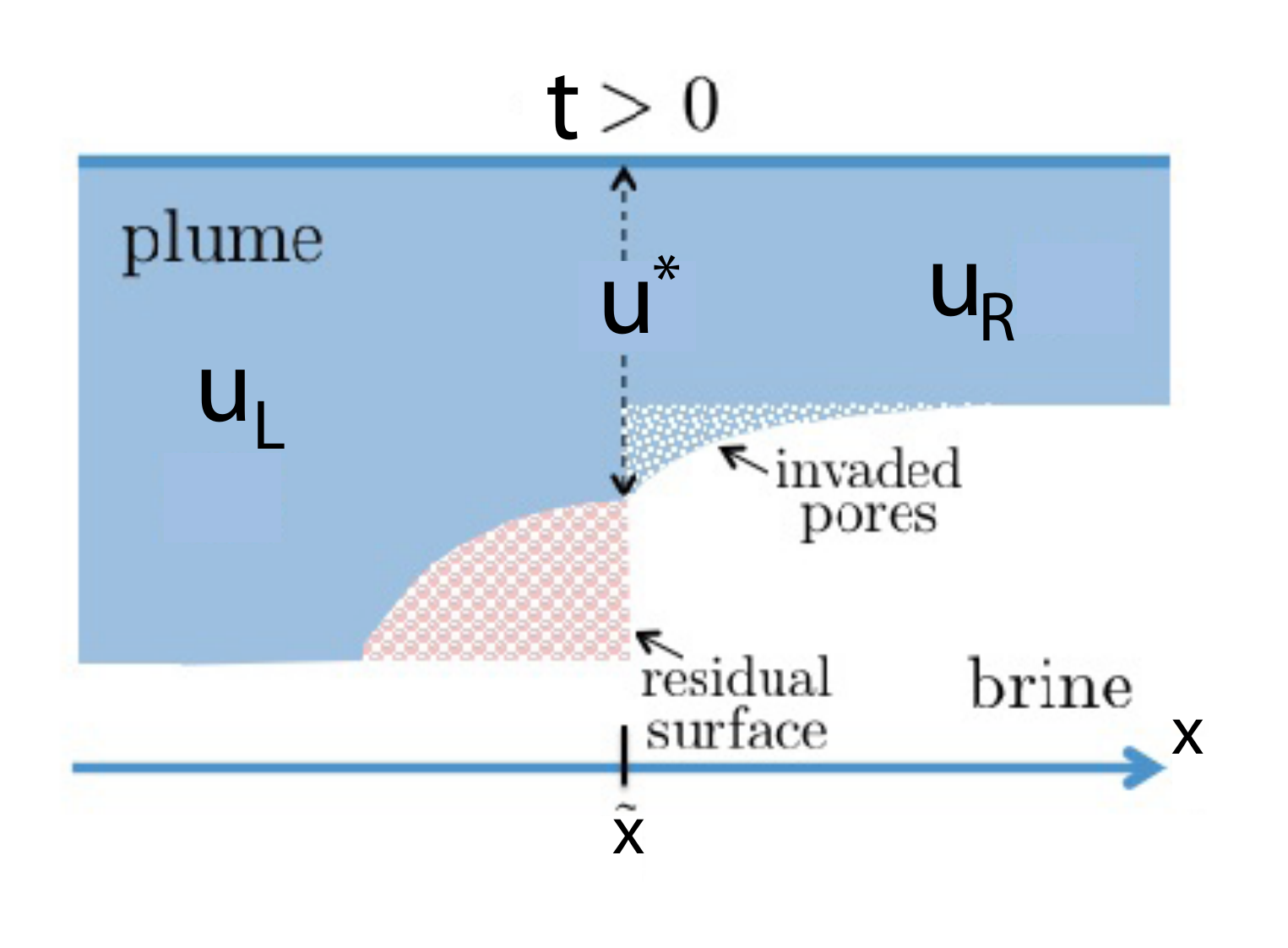} 
\\
 \hspace{1in} (a) \hspace{2.7in} (b) 
\end{tabular}
\caption{Rarefaction wave. (a) Left and right moving sections in the flux. (b)  CO$_2$ plume propagating left and right.}
\label{rare}
\end{figure}
The rarefaction solution \eqref{raresoln} varies continuously from $\ueta_L$ to $\ueta_R$ in Fig.\,\ref{exprare}(b). In particular, 
$\ueta(\xxi/\ttau)$ is continuous across $\xxi=0$ even though $\sigma$ in \eqref{raresoln} has a discontinuity at this position \cite{hesse08}.  Correspondingly, there is a discontinuity in the slope of the plume interface due to the jump in $\sigma\,.$ 
There is a jump $ [\ueta_\xxi]$ in the derivative $\partial_x\ueta$ at $\xxi=0,$ where $f'(\ueta)=0,$ and $\sigma$ switches from $\sigma=1$ to $\sigma=1-\epsilon.$ We calculate it assuming $f''(\ueta^*)<0:$ 
\begin{equation}
    [\ueta_\xxi]= \frac{1}{\ttau\epsilon f''(\ueta^*)}=-\frac{\sqrt{{\cal{M}}}}{2\ttau\epsilon}
\Label{etaxi}
\end{equation}
where the final equality uses the specific flux function \eqref{flux}. 

\section{Wave Interactions} \Label{waveinteractions}

\noindent We consider the Riemann problem, consisting of equation \eqref{PDE} with jump initial data
\begin{equation*}
\ueta(\xxi,0) =\left\{ \begin{array}{ll}
\ueta_L\,, & \quad \xxi<0\\[1pt]
\ueta_R\,, & \quad \xxi>0.
\end{array}\right.
\end{equation*}
It follows from \S\ref{charsshocks} that the solution
is an admissible shock if \,$\ueta_L<\ueta_R$ and a rarefaction fan 
if \,$\ueta_L>\ueta_R$. While the structure of these individual waves depends on the details of two flux functions and the switch between them, the outcome is, broadly speaking, the same as for a convex scalar conservation law with a single flux.

In this section, we consider pairs of Riemann problems. Each Riemann problem generates  a single wave; we are interested in whether the waves interact, and the result of the interaction. The results have significant differences from the corresponding wave interactions for a scalar equation with a single convex flux.

While a detailed classification is complicated, we focus on the main features of solutions of initial value problems with jump initial data of the form
\begin{equation}
\ueta(\xxi,0) =\left\{ \begin{array}{ll}
\ueta_L\,, & \quad \xxi<\xxi_1\\[3pt]
\ueta_M\,, & \quad \xxi_1<\xxi<\xxi_2\\[3pt]
\ueta_R\,, & \quad \xxi_2<\xxi\,
\end{array}\right.
\Label{LMR}
\end{equation}
in which $\ueta_L$ and $\ueta_R$ are different from $\ueta_M.$ 
Similar to the classical case, if $\ueta(\xxi,0)$ is decreasing, i.e. \textcolor{black}{$\ueta_L>\ueta_M>\ueta_R\,,$} then the solution consists of two rarefaction waves that do not approach. \textcolor{black}{Consequently, since the speed of an approximating expansion shock is between the speeds of the corresponding rarefaction's trailing and leading characteristics, two expansion shocks will not approach.} We treat the three remaining 
cases in turn\textcolor{black}{, and, if the data has an initial rarefaction, we examine the interactions involving expansion shock approximations. }


\vskip11pt 
\subsection{Case A: \textcolor{black}{Shock - Rarefaction:}  $\ueta_L<\ueta_M$ and $\ueta_R<\ueta_M$} \Label{CaseA}

In this case, we have a shock with speed $\Lambda$ emanating from $\xxi=\xxi_1$ at time $\ttau=0,$ and a rarefaction centered at $\xxi=\xxi_2>\xxi_1, \ttau=0.$ To see that the two waves approach, we check that the shock speed is greater than the speed 
of the trailing characteristic in the
rarefaction. There are two cases to consider. 
In case (i), $\Lambda>0\,,$ the shock admissibility condition requires $f'(\ueta_M)<\Lambda,$ so that 
the speed $\sigma f'(\ueta_M)$ of the trailing edge of the rarefaction is less than the shock speed, whether $\ueta_M<\ueta^*,$ for which $\sigma =1-\varepsilon,$ or $\ueta_M>\ueta^*,$ for which $\sigma =1.$ In case (ii), $\Lambda<0,$ so that $\ueta_\ttau>0$ and $\sigma=1-\epsilon.$ Thus, $\ueta_M>\ueta^*$ but now shock admissibility requires $\sigma f'(\ueta_M)<\Lambda,$ and the rarefaction, with trailing edge traveling at speed $f'(\ueta_M)<\sigma f'(\ueta_M)<\Lambda<0,$ approaches the shock.  


In Fig.~\ref{catalogclean2d}\textcolor{black}{,} we illustrate the solution as the interaction 
proceeds. In this and other figures, we plot exact solutions using the specific flux \eqref{flux} for illustration. On the left we show the track of the rarefaction through the flux curves as the characteristics fan from negative to positive speed.
The rarefaction fan provides the values of $\ueta$ on the right of the shock as the evolution proceeds.  The shock speed is represented by the slope of the chords in Fig. \ref{catalogclean2d}(a). As the speed switches from negative to positive, the chord moves from the lower flux graph to the upper, as $\ueta_\ttau$ changes sign. The crossover is represented by the horizontal dashed lines. In this example, the construction proceeds until the rarefaction wave has been completely absorbed by the shock. 
Since the initial data have $\ueta_R>\ueta_L,$ the long-time behavior is a single shock joining $\ueta_L$ to $\ueta_R.$ 
On the other hand, if $\ueta_L>\ueta_R\,,$ then the long-time behavior would be a rarefaction wave, the remnants of the short-time wave joining $\ueta_M$ to $\ueta_R,$ after the interaction with the shock wave has completed.  

This interaction of a shock with a rarefaction, illustrated in In Fig.~\ref{catalogclean2d}, appears to be similar to such interactions for a scalar conservation law with convex flux. However, there is a  significant difference. While the shock has negative speed, it is calculated from the flux $(1-\varepsilon)f(\ueta).$ \textcolor{black}{The shock} is admissible because the characteristics on the left have positive speed, and the faster characteristics on the right have speed $(1-\varepsilon)f'(\ueta),$ which is slower than the shock speed, as \textcolor{black}{shown in Fig. \ref{catalogclean2d}(b)}. In fact, for the smaller flux (in the lower graph), the shock satisfies the Lax entropy condition. However, as the shock turns and gains positive speed, we switch to the upper \textcolor{black}{flux curve}. The characteristics on the right both have negative speed to start with, and hence impinge on the shock. On the left, both characteristics travel faster than the shock. In fact, as the shock turns\textcolor{black}{,} it has zero speed, and the characteristics on the left for both fluxes have positive speed, so this property persists for some further time. 

However, as the shock continues to accelerate, there is a time, corresponding to shock location \textcolor{black}{, \raisebox{.5pt}{\textcircled{\raisebox{-.9pt} {1}}} in Fig. \ref{catalogclean2d}(b),} when \textcolor{black}{the shock} moves with the characteristic speed of $(1-\varepsilon)f'(\ueta_L)$ of the smaller flux, see the inclined dashed lines in Fig. \ref{catalogclean2d}(a) corresponding to $\ueta=\ueta_\text{graze}$. Consequently, if we continue to consider  only  the single slower family of characteristics (that were significant for the shock when it had negative speed), then the shock would fail to satisfy the Lax entropy condition at this time. By including the characteristics of the larger flux (which has already been invoked to calculate the shock speed) we retain admissibility of the shock. This device is consistent with causality, as the constant value of $\ueta$ is carried by both families of characteristics. 
This example and other similar instances are the reason for including both families of characteristics  {(hence, \textit{cross-hatch characteristics})} in open regions of the $(\xxi,\ttau)$ plane where $\ueta$ is constant.

\begin{figure}[!h]
\begin{center}
\begin{tabular}{cc}
{\includegraphics[height=3.2in]{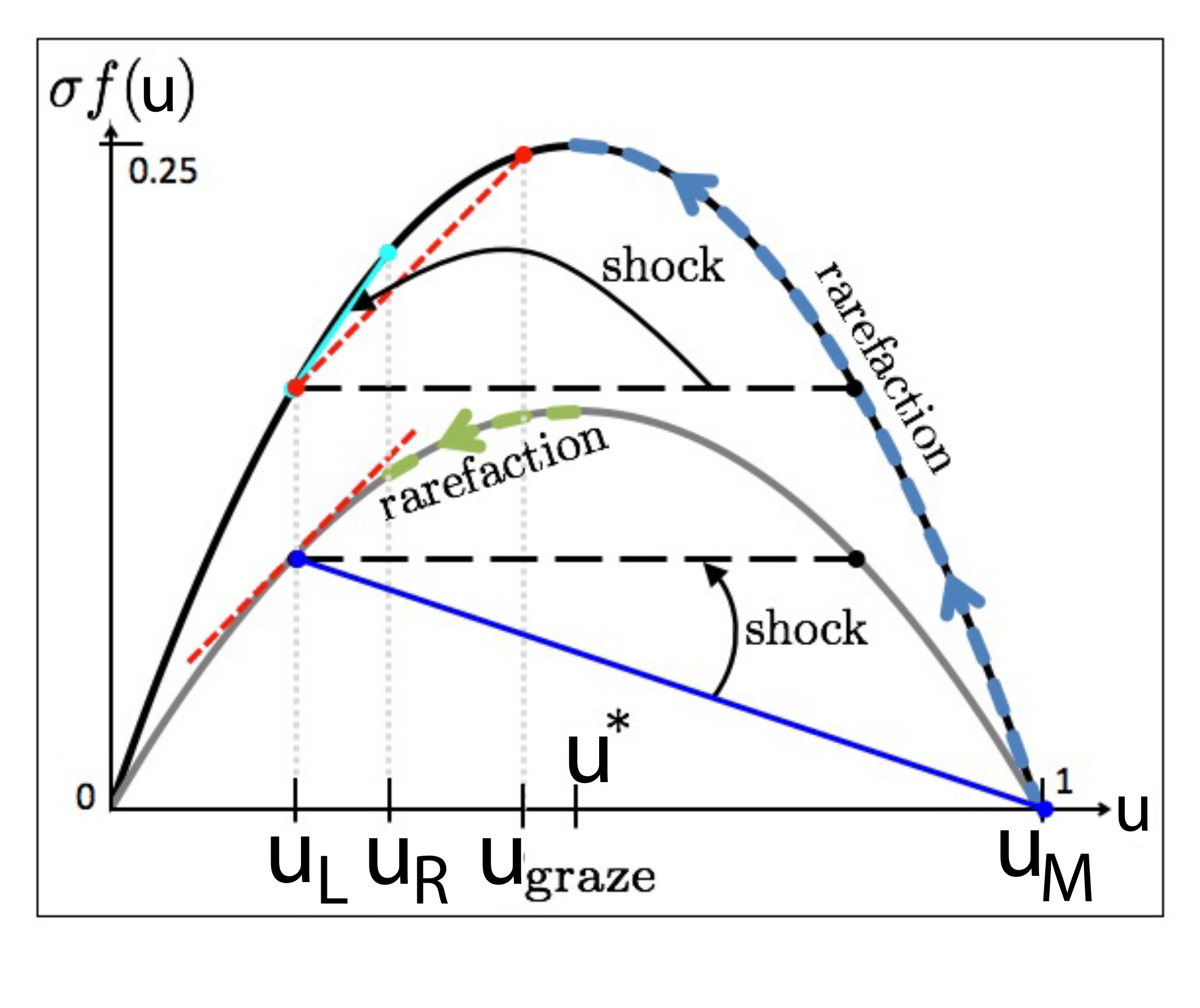}} & {\includegraphics[height=3.5in]{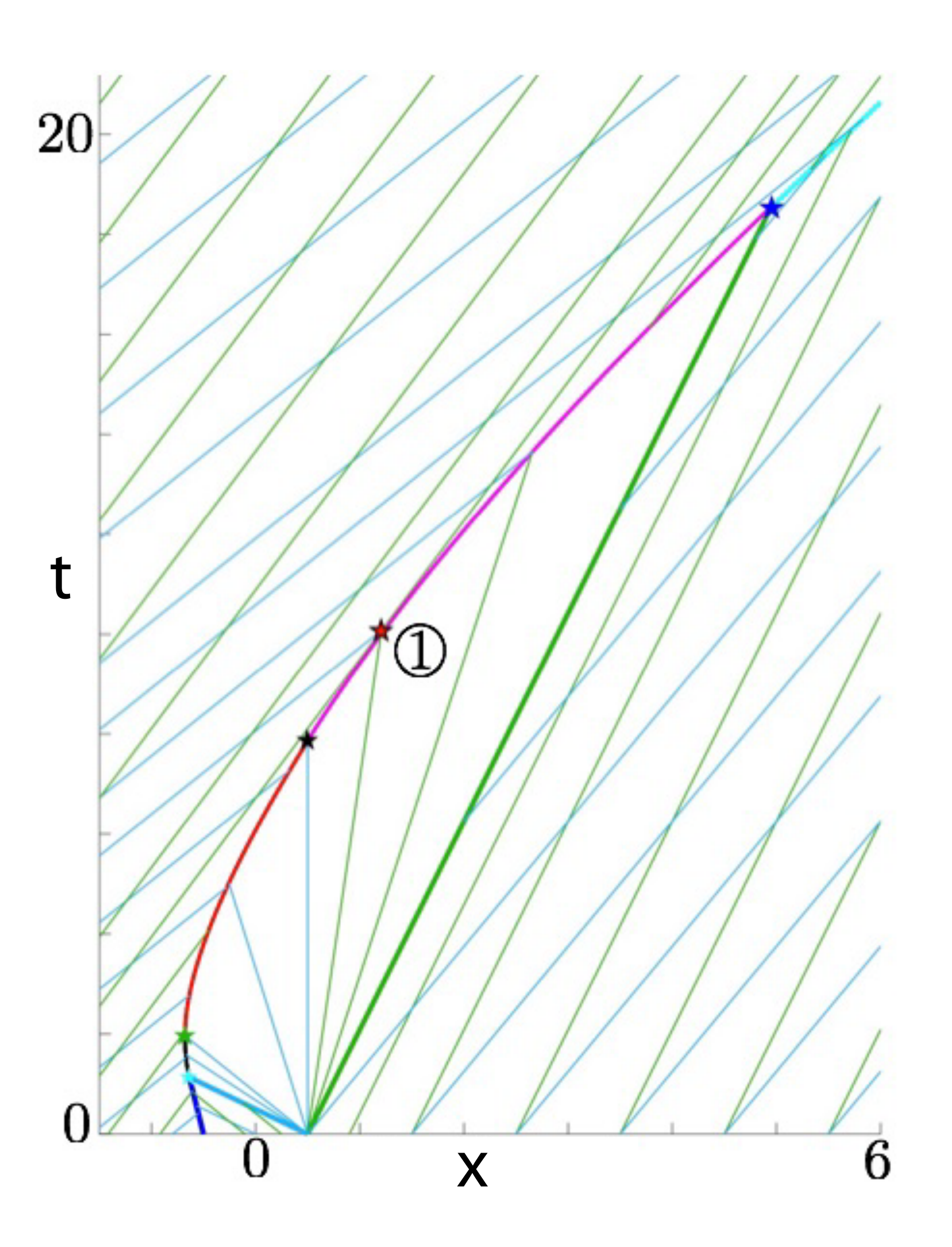}} \\[2pt]
\,\;\;(a) & \;\;\,\;\;(b) 
\end{tabular}
\caption{\textcolor{black}{Case A: Shock-rarefaction interaction with $\mathcal{M}=1,\,  \varepsilon=0.4, \, \ueta_L=0.2,\, \ueta_M=1,\, \ueta_R = 0.3.$ \ (a) Characteristic and shock speeds as the evolution proceeds. (b)  Characteristics and shock. Inclined dashed lines in (a) correspond to the  point \raisebox{.5pt}{\textcircled{\raisebox{-.9pt} {1}}} in (b) where a slower characteristic on the left grazes the shock.}}
\Label{catalogclean2d}
\end{center}
\end{figure}

\renewcommand\arraystretch{1.3} 

\vskip11pt 
\subsection{Case B: Shock - Shock: $\ueta_L<\ueta_M<\ueta_R$}
The second case involves a shock from a left state up to a middle state followed by a shock from the middle state up to a right state. Since the flux function is concave, the shock from $\ueta_L$ to $\ueta_M$ will have a greater shock speed than the shock from $\ueta_M$ to $\ueta_R$, so the shocks will approach each other and interact at a finite time to yield a single shock from $\ueta_L$ up to $\ueta_R$ with strength $\ueta_R-\ueta_L$. If the speeds of the approaching shocks have the same sign, the resulting shock has the same direction; if not, the resulting shock is forward if $f(\ueta_L)<f(\ueta_R)$\,, stationary if $f(\ueta_L)=f(\ueta_R)$\,, or backward if $f(\ueta_L)>f(\ueta_R)$\,. The total variation is unchanged before and after the discontinuities interact, and the middle state is eliminated in finite time. 



\vskip11pt 
\subsection{Case C: Rarefaction - Shock: $\ueta_M<\ueta_L$ and $\ueta_M<\ueta_R$}
This case mirrors Case A, in that the short-time solution is a rarefaction wave to the left of a shock wave. However, whereas in Case A the two waves approach, in Case C their approach depends on further restrictions on the data. The reason for this is that the slower characteristics on the left can leave the shock (Lemma~\ref{lemma1}); they are necessarily parallel to the leading edge of the rarefaction.  We distinguish two subcases in which the waves do not approach: 

(i) If $\ueta_M\textcolor{black}{ \leq }\ueta^*,$ define 
$\tilde\ueta_M$ by 
$$
\frac{f(\tilde\ueta_M)-f(\ueta_M)}{\tilde\ueta_M-\ueta_M}=(1-\varepsilon)f'(\ueta_M),
$$
shown in Fig. \ref{Cnoint}, and let $\lambda_M$ denote this speed. Then $\lambda_M>0$ is the speed of the leading edge of the rarefaction, and if $\ueta_R=\tilde\ueta_M,$ then it is also the speed of the shock, since the shock has a jump up and positive speed. Then for $\ueta_M\textcolor{black}{\leq}\ueta^*$ and  
\beq\label{noint1}
\ueta_M<\ueta_R\leq\tilde\ueta_M, \quad \ueta_M<\ueta_L,
\eeq
the shock from $\ueta_M$ to $\ueta_R$ has positive and larger speed:
$$
\frac{f(\ueta_R)-f(\ueta_M)}{\ueta_R-\ueta_M}\geq\lambda_M.
$$
Thus, \eqref{noint1} is sufficient to guarantee that the shock and rarefaction do not approach.

\begin{figure}[!t]
\begin{center}
\begin{tabular}{cc}
{\includegraphics[height=2.5in]{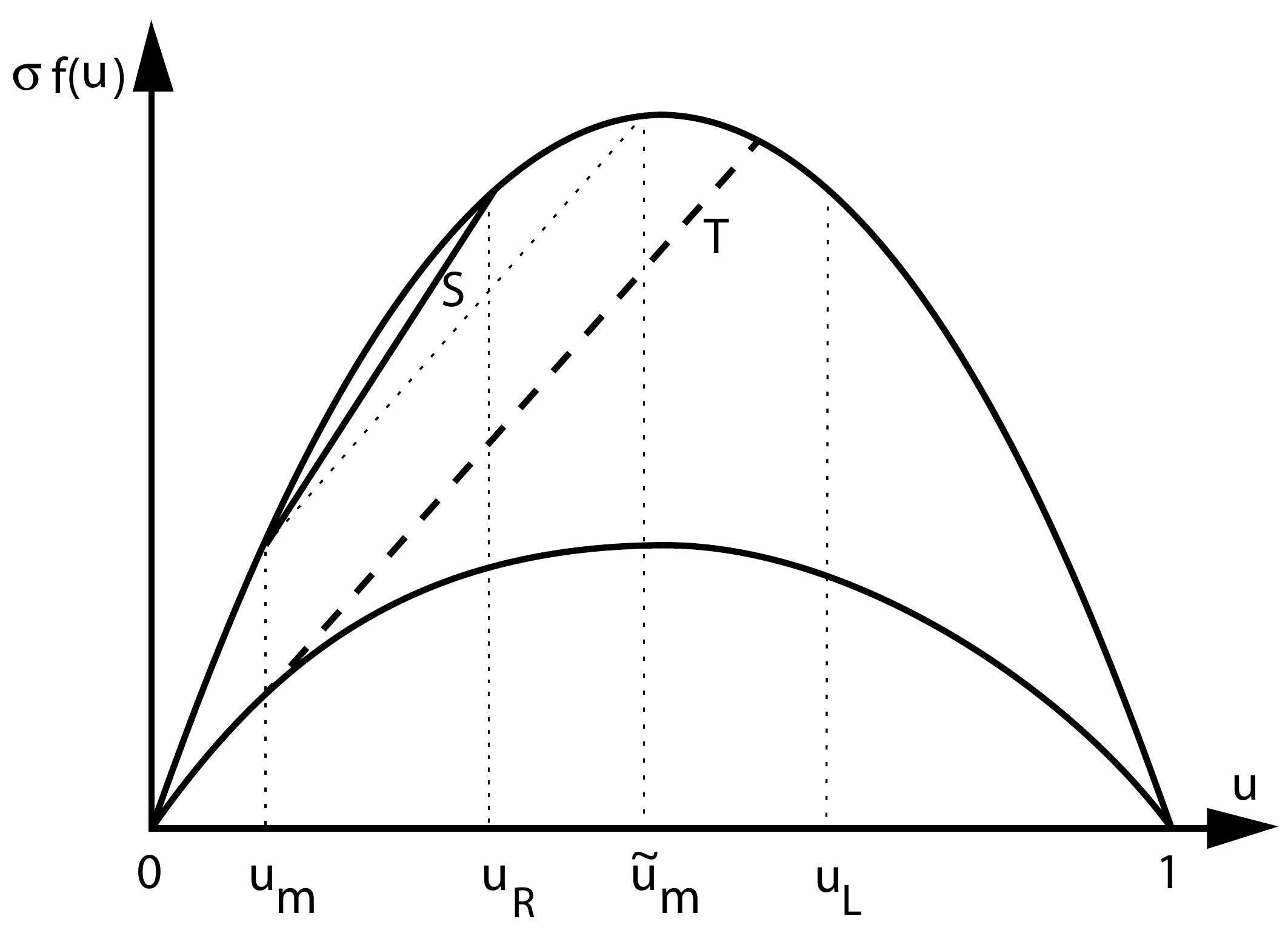}} &{\includegraphics[height=2.2in]{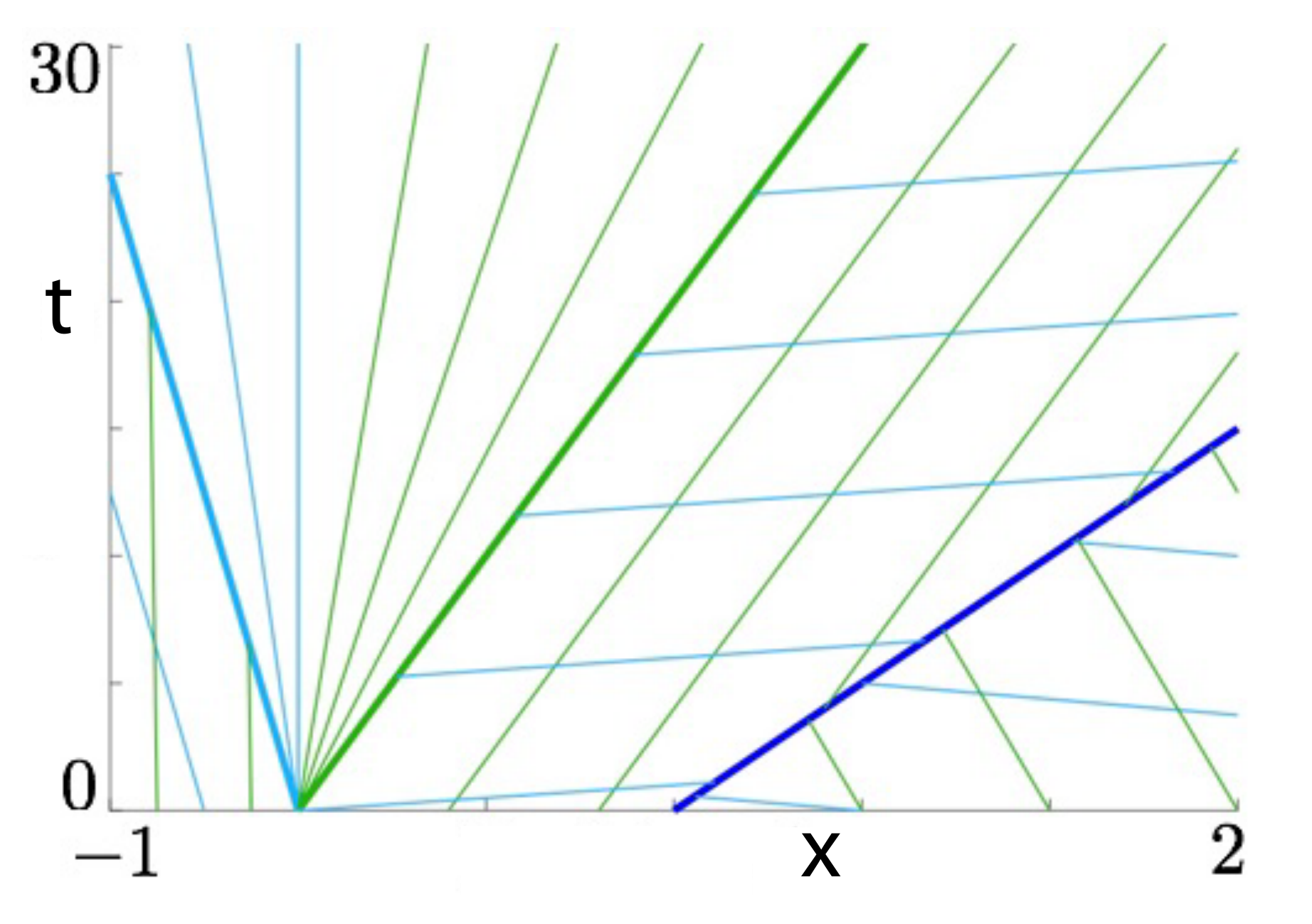}} \\[2pt]
\,\;\;(a) &\quad \,\;\;(b) 
\end{tabular}
\caption{Case C(i): Rarefaction and shock do not interact. (a) shock speed: slope of solid chord S; rarefaction leading edge speed: slope of dashed tangent T.   (b) $\xxi,\ttau$ plane with characteristics.
}
\Label{Cnoint}
\end{center}
\end{figure}

\vskip11pt
(ii) Similarly, if $\ueta_M>\ueta^*,$ then the shock speed and speed of the leading edge of the rarefaction wave are both negative. In this case, the rarefaction is backward and uses the larger flux $f(\ueta)$ whereas the shock uses the lower flux $(1-\varepsilon)\textcolor{black}{f(\ueta)}\,.$ \textcolor{black}{C}onsequently, the \textcolor{black}{interaction} condition becomes 
$$
(1-\varepsilon)\,\frac{f(\ueta_R)-f(\ueta_M)}{\ueta_R-\ueta_M}>f'(\ueta_M)\,.
$$
Define $\overline\ueta_M>\ueta_M$ by
$$
\left\{\barr{l}
\overline\ueta_M=1 \,, \ \ \ \mbox{if}\quad (1-\varepsilon)\,\ds\frac{f(1)-f(\ueta_M)}{1-\ueta_M}>f'(\ueta_M)\\[11pt]
(1-\varepsilon)\,\ds\frac{f(\overline{\ueta}_M)-f(\ueta_M)}{\overline{\ueta}_M-\ueta_M}=f'(\ueta_M)\,,\ \ \mbox{otherwise.}
 \earr
 \right.
 $$
Then the two waves do not approach if $\ueta^*<\ueta_M$ and 
\beq\label{noint2}
\ueta_M<\ueta_R\leq \overline\ueta_M\,, \quad \ueta_M<\ueta_L.
\eeq
In summary, if neither \eqref{noint1} nor \eqref{noint2} are satisfied by $\ueta_R\,,$ then the rarefaction wave and shock wave interact much as in Case~A, see Fig. \ref{catalogclean6e}(a). Otherwise, the shock travels faster than the rarefaction, and there is no interaction, as in Fig. \ref{Cnoint}(b).

\textcolor{black}Unlike Cases A and B, not all initial conditions in Case C lead to an eliminated initial middle state in finite time. Some solutions in Case C exhibit unusual behavior, due to the flux discontinuity, that does not arise in scalar equations with a single flux: shock speeds determined by one flux curve can equal corresponding characteristic speeds found on the other flux curve. In Fig. \ref{catalogclean6e}(b), the plume asymptotically approaches a height of $\tilde{\ueta}\in\Big(\max(\ueta^*,\ueta_M)\,,\,\min(\ueta_L,\ueta_R)\Big)$ such that 
\begin{equation*}
f'(\tilde{\ueta})=\sigma \,\dfrac{f(\ueta_R)-f(\tilde{\ueta})}{\ueta_R-\tilde{\ueta}}.
\label{etatilde}
\end{equation*}
Hence, if $\ueta_L\geq\tilde{\ueta}$, the backward shock does not reach the rarefaction's trailing characteristic; the shock speed approaches the characteristic speed corresponding to $\tilde{\ueta}$ 
labeled e in Fig. \ref{catalogclean6e}(b). The result approaches a rarefaction from $\ueta_L$ down to $\tilde{\ueta}$ and a shock from $\tilde{\ueta}$ up to $\ueta_R$; since $\ueta_M<\tilde{\ueta}$, the total variation of the solution decreases to $\ueta_L+\ueta_R-2\,\tilde{\ueta}$. 

However, if $\ueta_L<\tilde{\ueta}$ as in Fig. \ref{catalogclean6e}(a), the middle state is eliminated in finite time, resulting in a decrease of total variation to $\ueta_R-\ueta_L$. It is also possible for a middle state to asymptote to a value $\overline{\ueta}\in\Big(\ueta_M\,,\,\min(\ueta^*,\ueta_L,\ueta_R)\Big)$ such that 
\begin{equation*}
\sigma\,f'(\overline{\ueta})=\dfrac{f(\ueta_R)-f(\overline{\ueta})}{\ueta_R-\overline{\ueta}}\,
\label{etabar}
\end{equation*}
since the speed of a forward shock is determined by the upper flux curve, and the characteristic speed to the right of the center of a rarefaction is found on the lower flux curve. Again, the total variation of the solution decreases. Hence, for Case C, if there is an interaction, the total variation always decreases, and the number of outgoing waves is non-increasing. 

\renewcommand\arraystretch{0.7}
\begin{figure}[!h]
\begin{center}
\begin{tabular}{c}
{\includegraphics[height=2.3in]{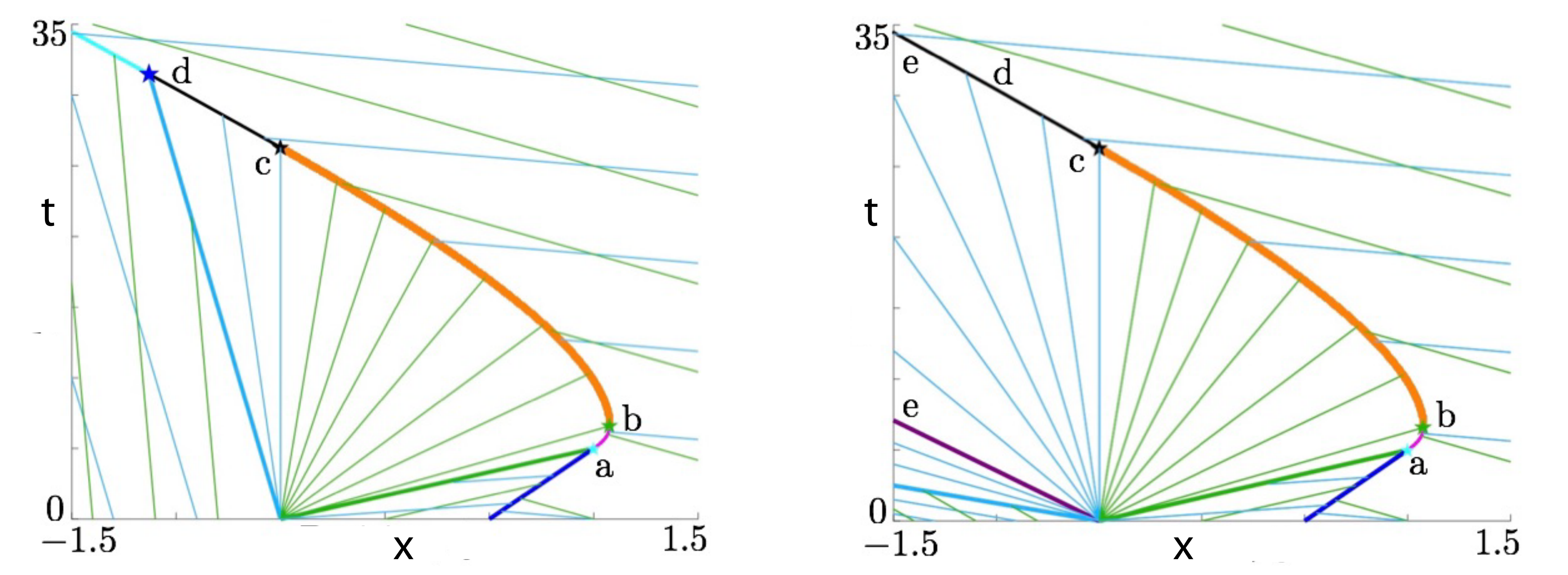}} \\
(a) \hskip3.2in (b)
\end{tabular}
\caption{\textcolor{black}{Case C(ii): Rarefaction-shock interactions, with   $\mathcal{M}=1, \, \varepsilon = 0.7, \, \ueta_M=0,\,\ueta_R=0.9. $ \  (a)  $\ueta_L=0.51<\tilde{\ueta}$\,, and (b) $\ueta_L=0.7>\tilde{\ueta}$\,.}}
\Label{catalogclean6e}
\end{center}
\end{figure}
\subsection{Interactions of Shocks and Expansion Shocks}

\subsubsection{\textcolor{black}{Case a: Shock - Expansion Shock: $\ueta_L<\ueta_M$ and $\ueta_R<\ueta_M$}}

 This sub-case corresponds to case A above, in which a shock necessarily interacts with a rarefaction wave. 
 However, when the rarefaction is replaced by a piecewise constant approximation consisting of expansion shocks, the shock from $\ueta_L$ up to $\ueta_M$ may not meet the slowest expansion shock on the right. The two waves move apart if the shock has positive speed and the expansion shock has larger speed, or if the shock has negative speed and the expansion shock has either less negative or positive speed. To analyze the situation, we consider the expansion shock from $\ueta_M$ to $\ueta_R.$ 
  Let $\Lambda_{LM} = \dfrac{f(\ueta_M)-f(\ueta_L)}{\ueta_M-\ueta_L}.$}
\vskip5pt
 (i) If $\Lambda_{LM} >0$\,, a forward shock with speed $\Lambda_{LM}$ connects $\ueta_L$ and  $\ueta_M$\,. If $\ueta_M\leq \ueta^*$, the expansion shock between $\ueta_M$ and $\ueta_\RR$ also has positive speed; however, if $\ueta_M>\ueta^*$\,, the expansion shock could have positive, zero, or negative speed. The shock and expansion shock will move apart only if the expansion shock moves faster than the  shock, in which case,  the expansion shock has speed
  $\Lambda_{MR}= (1-\epsilon)\dfrac{f(\ueta_M)-f(\ueta_R)}{\ueta_M-\ueta_R}.$ Let $\tilde{\ueta}_M<\ueta_M$ be such that
\begin{equation*}
\left\{\barr{l}
\tilde{\ueta}_M=0\,, \ \ \ \mbox{if}\quad   (1-\varepsilon)\,\ds\frac{f(\ueta_M)-f(0)}{\ueta_M}< \Lambda_{LM}\\[11pt] 
(1-\varepsilon)\,\ds\frac{f({\ueta}_M)-f(\tilde{\ueta}_M)}{{\ueta}_M-\tilde{\ueta}_M}=\Lambda_{LM}\,,\ \ \mbox{otherwise.} 
 \earr
 \right.
\end{equation*}
Hence, the shock and approximating expansion shock(s) do not approach if $\ueta_L<\ueta_M$ with $\Lambda_{LM}>0$\, and $\ueta_\RR$ satisfies 
\begin{equation} \Label{shockexp1}
0\leq {\ueta_\RR\leq \tilde{\ueta}_M}<\ueta_M \,. 
\end{equation}

\textcolor{black}(ii) For the case when $\Lambda_{LM}\leq 0$, the shock speed $(1-\varepsilon)\Lambda_{LM}$ is reduced due to residual trapping. Define $\overline\ueta_M<\ueta_M$ to be such that
\begin{equation*}
(1-\varepsilon)\Lambda_{LM} = \dfrac{f(\ueta_M)-f(\overline\ueta_M)}{\ueta_M-\overline\ueta_M}\,.
\end{equation*}
The shock wave and expansion shock wave do not interact if $\ueta_L<\ueta_M$ has 
$\Lambda_{LM}\leq 0$ and $\ueta_\RR$ is such that
\begin{equation}\Label{shockexp2}
0\leq \ueta_\RR\leq \overline{\ueta}_M<\ueta_M .
\end{equation}
Hence, if we have an expansion shock from $\ueta_M$ down to $\ueta_\RR$, where $\ueta_\RR$ does not satisfy \eqref{shockexp1} or \eqref{shockexp2}, then the shock and   expansion shock collide, producing a single admissible shock from $u_L$ to $u_R.$ 
 
\subsubsection{\textcolor{black}{Case c: Expansion Shock - Shock: $\ueta_M<\ueta_L$ and $\ueta_M<\ueta_R$}}

This case is analyzed similarly to Case b, with conditions for the approach or separation of the two waves, analogous to Case C, where a rarefaction wave to the left of a shock may fail to approach the shock because the fastest characteristic in the rarefaction is slower than the shock speed. Correspondingly, when a rarefaction from $u_L$ to $u_M<u_L$  is approximated with one (or more) expansion shock(s), the fastest (right-most) expansion shock connects $\ueta_E\in (\ueta_M\,,\,\ueta_L\,]$ down to $\ueta_M$\, with speed 
$\Lambda_{EM} =\sigma_L \dfrac{f(\ueta_E)-f(\ueta_M)}{\ueta_E-\ueta_M}.$ If this speed is  less than 
$\Lambda_{RM} = \sigma_R\dfrac{f(\ueta_R)-f(\ueta_M)}{\ueta_R-\ueta_M},$ then the two waves fail to interact and all the expansion shocks approximating the rarefaction move away from the shock.  Here, $\sigma_L=1$ if and only if $\Lambda_{EM}<0,$ and $\sigma_R=1$ if and only if $\Lambda_{RM} > 0.$

\textcolor{black} 
If both waves are moving right, then they   interact only if $u_R>u_L.$ If they do interact, then the result is an admissible shock from $u_L$ to $u_R$. A similar argument applies to left-moving waves: either they separate, or the result is an admissible shock from $u_L$ to $u_R$. Consequently, the number of waves either remains at two, with no change in the total variation, or is decreased to one, with a corresponding decrease in total variation.



The overall result of binary interactions between shock waves and expansion shocks  is that the  total variation and number of waves decreases, but adjacent waves may move apart.

%
%

\section{Initial Value Problems} \Label{existenceproof}
Plume migration within a porous aquifer depends on the geometry of the carbon dioxide plume at the end of injection \cite{huppert14}, \cite{juanes10}. An analytic solution for a specific idealized CO$_2$ plume is constructed by Hesse, Orr, and Tchelepi \cite{hesse08}. In this section
we consider the scalar conservation law \eqref{PDE} with a general initial plume of supercritical carbon dioxide,
\begin{align}
\left\{ \begin{array}{ll}
\ueta_\ttau+\big(\sigma\,f(\ueta)\big)_\xxi = 0\,, & \quad\xxi\in\mathbb{R}\,,\; \ttau>0\,,\\[6pt]
\ueta(\xxi,0)=\ueta_0(\xxi)\,, & \quad \xxi\in\mathbb{R}\,,
\end{array}\right.
\Label{cauchy}
\end{align}
in which $\ueta_0 \in L^1(\mathbb{R})\, \cap\, BV(\mathbb{R})\,$ with $\,0\leq  \ueta_0 \leq1$\,.

\subsection{Wave-Front Tracking}

Dafermos \cite{dafermos72} introduced wave-front tracking as a method to construct approximate solutions for scalar, nonlinear partial differential equations. The method has since been greatly generalized to systems of hyperbolic conservation laws \cite{bressan91}, \cite{bressan00}.
 In this section, we describe wave-front tracking, following the approach of LeFloch \cite{lefloch}.
In the wave-front tracking algorithm, we first approximate the initial plume shape with a sequence of piecewise constant functions, $\ueta_0^h(\xxi), \ h>0$, such that 
\begin{align} \Label{initialapprox}
&\inf(\ueta_0) \leq \ueta_0^h \leq \sup(\ueta_0) \nonumber\\
&TV(\ueta_0^h)\leq TV(\ueta_0)\\ 
&\ueta_0^h\rightarrow \ueta_0 \text{ in  }L^1 \text{ as } h\rightarrow0^+. \nonumber 
\end{align} 
  Each approximation $\ueta_0^h$ is constructed to have a finite number of discontinuities. 
  The construction of a piecewise-constant solution for short time involves solving the Riemann problems associated with each discontinuity in $\ueta_0^h.$ Rarefaction waves are replaced by a finite number of expansion shocks of magnitude $h$ or less. When waves meet, we refer to the  collision as an interaction. Each interaction results in a   Riemann problem in which the initial jump may exceed the threshold $h.$   If the resulting solution is an admissible shock, it is propagated forward without change. If the resulting solution is a rarefaction wave, then (as observed in the previous section) the magnitude is necessarily smaller than $h;$ it is approximated by an expansion shock, traveling with the shock speed of that discontinuity.  Continuing in this way, we generate  a piecewise constant solution of the conservation law.



In \S\ref{waveinteractions}, we showed that the number of waves and total variation decreased or remained constant at any interaction. Consequently, since there are  finitely many   discontinuities initially, there are a finite number of interactions and  no accumulation points. 
Thus, the number of wave interactions and resulting wave-fronts in each $\ueta^h$ remains finite for all $\ttau>0$, so the approximations are well defined globally in time. 
\cite{bressan00}.  
  
 As observed in the previous section, the total variation is non-increasing, and each approximation $\ueta^h(\xxi,\ttau)$ is bounded by $\ueta_0^h(\xxi)$\,. It follows from \eqref{initialapprox} that, at any position and time, \,$\inf(\ueta_0)\leq \ueta^h(\xxi,\ttau)\leq \sup(\ueta_0)$\,; hence,\, $\left\|\ueta^h(\xxi,\ttau)\right\|_{L^\infty} \leq 1$\,. We also have \,$TV\big(\ueta^h(\cdot,\ttau)\big) \leq TV\big(\ueta_0(\cdot)\big)$\, for all $\ttau>0$\,.

\vskip11pt
Since we have established that there are a finite number of waves, there will be a finite number, $k$, of classical and expansion shocks in $\ueta^h$ within $\big[\,\ttau_1,\ttau_2\,\big]$, any time interval containing no interaction time. For $m=1, \ldots, k$\,, let $y_m'$ be the speed of propagating shock front $\xxi=y_m(\ttau)$ in $\ueta^h$ for $\ttau\in\big[\,\ttau_1,\ttau_2\,\big]$\,; by \eqref{flux}, $\left|y_m'\right|\leq \sup\left|f' \right|<\infty$\,. The approximate solution to the left/right of wave-front $ y_m$ is $\ueta^h\big(y_m(\ttau)^\mp,\ttau\big)$. We following estimate is based on the change in area under the graph of $\ueta^h(\ttau)$ due to the motion of individual waves.  
\begin{align*}
\Big\|\ueta^h(\xxi,\ttau_2)-\ueta^h(\xxi,\ttau_1)\Big\|_{L^1} & \;\leq \; \sum_{m=1}^k \left|\ueta^h\big(y_m(\ttau_1)^-,\ttau_1\big)-\ueta^h\big(y_m(\ttau_1)^+,\ttau_1\big)\right| \left|y_m' \right| \big| \ttau_2-\ttau_1\big|\\[2pt]
& \;\leq \;TV(\ueta_0) \sup\left|f' \right| \big| \ttau_2-\ttau_1\big|.
\end{align*}
We have shown that conditions for both  Helly's Theorem and the time-dependent version (\cite{bressan00},\cite{lefloch}) are satisfied.   Hence, there exists a subsequence of $\ueta^h$, which we also label $\ueta^{h}$, and a BV function   $\ueta: \mathbb{R}\times \mathbb{R}^+\rightarrow [\,0\,,1\,]\,$, such that
\begin{align} \Label{helly}
&\ueta^h(\xxi,\ttau)\rightarrow \ueta(\xxi,\ttau) \quad\qquad \text{in } L_\text{loc}^1\,, \nonumber\\[6pt]
&\big\|\ueta(\ttau)\big\|_{L^\infty}+TV\big(\ueta(\ttau)\big) \leq \kappa\,, \text{ and}\\[6pt]
&\big\|\ueta(\ttau_2)-\ueta(\ttau_1)\big\|_{L^1}\leq \kappa\,\big|\ttau_2-\ttau_1\big|,\nonumber 
\end{align}
for all $\xxi\in\mathbb{R}$\,, $\ttau,\,\ttau_1,\,\ttau_2\in\mathbb{R}^+$, and some $\kappa>0$\,.


\vskip11pt
Combining \eqref{helly} with the lower semi-continuity property  $TV\big(\ueta(\cdot,\ttau)\big) \leq \displaystyle \mathop{{\lim\inf}}_{h \to 0^+} \,TV\big(\ueta^h(\cdot,\ttau)\big)$ we have  $TV\big(\ueta(\cdot,\ttau)\big) \leq  TV \big(\ueta_0(\cdot)\big)$ for all   $\ttau\geq0$\,. Similarly, since $\ueta^h$ 
converges to $\ueta$\,, it follows that 
$\,\inf(\ueta_0)\leq \ueta(\xxi,\ttau)\leq \sup(\ueta_0)$\,. We also have \,$\left[\ueta^h(\xxi,\ttau_2) - \ueta^h(\xxi,\ttau_1)\right]\rightarrow \big[\ueta(\xxi,\ttau_2) - \ueta(\xxi,\ttau_1)\big]$ \,in $L^1_\text{loc}$\, by \eqref{helly}, and it follows from the lower semi-continuity property of norms that \,$\big\|\ueta(\xxi,\ttau_2)-\ueta(\xxi,\ttau_1)\big\|_{L^1} \leq\mathop{{\lim\inf}}_{h \to 0^+}\big\|\ueta^h(\xxi,\ttau_2)-\ueta^h(\xxi,\ttau_1)\big\|_{L^1}\,.$ Finally, from the uniform estimate above, we have $\,\big\|\ueta(\xxi,\ttau_2)-\ueta(\xxi,\ttau_1)\big\|_{L^1} \leq TV(\ueta_0) \sup\left|f' \right| \big| \ttau_2-\ttau_1\big|\,$ for all $\ttau_1,\ttau_2\geq0$\,.

\vskip11pt

The wave-front tracking approximations $\ueta^h$ are exact solutions of \,$\ueta^h_\ttau+\big(\sigma\,f(\ueta^h)\big)_\xxi=0$ \,since  the Rankine-Hugoniot jump condition is satisfied across all classical and expansion shocks. However, we are unable to take the limit as $h\to 0+$ for a pair of reasons: First, we do not have a weak formulation of the Cauchy problem, and second, as $h\to 0,$ the value of $\sigma(\ueta^h_t)$ changes and it is not clear how to formulate the limit $\lim_{h\to0}\sigma(\ueta^h_t)f(\ueta^h)_\xxi$ in the sense of distributions, which should be $\sigma(\ueta_t)f(\ueta)_\xxi.$ If such problems can be resolved, then establishing that the limit is an entropy solution in the appropriate sense generalized to  the model is straightforward. 

\section{Discussion}
The   tracking a plume of supercritical carbon dioxide after it has been injected into a deep saline aquifer is modeled by a scalar partial differential equation that has unusual features due to the property of deposition of $CO_2$ bubbles as the plume migrates. In the model of Hesse et al \cite{hesse08}, this is achieved by reducing the flux by a constant scale as the plume migrates away from a region of space, leaving behind bubbles of sequestered $CO_2.$ In this paper, we have explored some interesting properties of the model that fall outside the conventional theory of conservation laws. 

The method of characteristics has an interesting twist, due to the presence of two characteristic speeds. Since the switch occurs when either $u_x$ or $f'(u)$  changes sign, tracking maxima and minima of   $u(x,t)$  the solution propagates either as a corner, or as an expanding interval in $x$ over which $u(x,t)$ is constant. Similarly, if a rarefaction wave includes values of $u$ that cross $u=u^*,$ where $f(u)$ has a maximum, then  the rarefaction wave includes a corner, where the slope $u_x$ jumps as $u$ crosses $u^*.$ 

In order to define shock waves, we have to generalize the Lax entropy condition in that the admissible behavior of characteristics on either side of the discontinuity has to be interpreted appropriately. A consequence is that each value of $u$ has two possible characteristic speeds, namely $f'(u)$ and $(1-\epsilon)f'(u).$  The choice depends on the direction of propagation of the wave, so that the choice switches if a shock wave changes direction. 
 To accommodate this behavior, we express admissibility in terms of both families of characteristics. 
 
 These phenomena associated with characteristics and shock waves appear when describing the interaction of pairs of waves. We find that shock-to-rarefaction interactions can be complete in finite time, leaving a shock wave, or can persist, resulting in a remaining rarefaction and a shock whose speed approaches characteristic speed. 
 
 The asymptotic behavior as $t\to \infty$ shown in Fig.~\ref{catalogclean6e}(b) suggests an unusual rarefaction-shock construction, in which $u=u_L$ is connected to $\tilde{u}<u_L$ by a rarefaction wave, whose fastest characteristic speed (the speed of the right-most characteristic in the rarefaction fan) is the same as the shock speed of a jump from $\tilde{u}$ to $u_R>\tilde{u}.$ This composite wave, in which the shock is characteristic on one side and does not decay, is unusual, because the flux functions are convex, whereas shock-rarefactions are expected to appear only when genuine nonlinearity fails; that is, for non-convex flux functions. 
 
 It would be interesting to know how the notion of weak solution can be formulated for equation \eqref{PDE}. Although it is clear how to treat piecewise smooth solutions, the convergence result from wave front tracking does not guarantee that the limit is a piecewise smooth function, even if the initial data are smooth. In terms of the application, it would be interesting to know whether compactly supported initial data collapses to zero in finite time, signifying the desirable property of complete sequestration in a finite time and over a finite distance. Of particular significance would be an estimate of the  maximum time over which this would occur, and the corresponding maximum distance any given plume would migrate before giving up all its $CO_2$ to sequestered bubbles.


\end{document}